\documentclass[11pt, oneside]{amsart}   	
\usepackage{geometry}
\usepackage{amsmath}            		
\usepackage{graphicx}				
\usepackage{amssymb}
\usepackage[all]{xypic}
\usepackage{bbm}  
\usepackage{slashed}

\newcommand{\Tr}{{\rm Tr}}
\newcommand{\dsl}{{\del\kern-10pt /}}
\newcommand{\Asl}{{A\kern-15pt /}}

\newcommand{\CG}{{\hbox{{$\mathcal G$}}}}
\newcommand{\CC}{{\hbox{{$\mathcal C$}}}}

\newcommand{\CE}{{\hbox{{$\mathcal E$}}}}

\newcommand{\CS}{{\hbox{{{$\mathcal S$}}}}}

\newcommand{\R}{\mathbb{R}}
\newcommand{\C}{\mathbb{C}}
\newcommand{\Z}{\mathbb{Z}}
\newcommand{\N}{\mathbb{N}}

\newcommand{\cX}{\mathfrak{X}}

\newcommand{\lbiprod}{{>\!\!\!\triangleleft\kern-.53em\cdot\,}}  
\newcommand{\rbiprod}{{\,\cdot\kern-.53em\triangleright\!\!\!<}}  

\newcommand{\note}[1]{}
\newcommand{\Hom}{{\rm Hom}}

\newcommand{\Mor}{{\rm Mor}}

\newcommand{\coev}{{\rm coev}}

\newcommand{\del}{\partial}
\newcommand{\ev}{{\rm ev}}

\newcommand{\extd}{{\rm d}}
\newcommand{\isom}{{\cong}}
\newcommand{\eps}{{\epsilon}}
\newcommand{\tens}{\mathop{\otimes}}
\newcommand{\la}{{\triangleright}}

\newcommand{\id}{{\rm id}}

\newcommand{\trace}{{\rm Tr}}

\newcommand{\dirac}{{ \slashed{D} }}

\newcommand{\pd}{\partial}

\newcommand{\rdy}{{\Lambda\kern-4.3pt\Lambda}}
\newcommand{\coH}{\mathrm{H}}
\renewcommand{\imath}{\mathrm{i}}
\newcommand{\kk}{\mathbbm{k}}

\newcommand{\und}[1]{{\underline {#1}}}

\newcommand*{\doublenabla}{%
  \nabla\mkern-12mu\nabla
}   

\newcommand{\eproof}{\hfill $\square$

\medskip}

\renewcommand\arraystretch{1.7}

\newcounter{piccie}
\begin{document}

\title[Quantum Bianchi identities]{Quantum Bianchi identities and characteristic classes via DG categories}

\author{Edwin J.\ Beggs \&\ Shahn Majid}
\address{EJB: Dept of Mathematics, Swansea University\\ Singleton Parc, Swansea SA2 8PP\\
SM: School of Mathematical Sciences, Queen Mary, University of London\\
 Mile End Rd, London E1 4NS}

\keywords{Quantum Riemannian geometry, noncommutative geometry, Chern-Connes pairing, Bianchi identity, tensor category, quantum group, q-sphere, quantum gravity, quantum spacetime}

\maketitle 

\newtheorem{theorem}{Theorem}[section]
\newtheorem{lemma}[theorem]{Lemma}
\newtheorem{proposition}[theorem]{Proposition}
\newtheorem{corollary}[theorem]{Corollary}
\newtheorem{definition}[theorem]{Definition}
\newtheorem{example}[theorem]{Example}
\newtheorem{conjecture}[theorem]{Conjecture}

\begin{abstract} We show how DG categories arise naturally in noncommutative differential geometry and use them to derive noncommutative analogues of the Bianchi identities for the curvature of a connection. We also give a derivation of formulae for characteristic classes in noncommutative geometry following Chern's original derivation, rather than using cyclic cohomology. We show that a related DG category for extendable bimodule connections is a monoidal tensor category and in the metric compatible case give an analogue of a classical antisymmetry of the Riemann tensor. The monoidal structure implies the existence of a cup product on noncommutative sheaf cohomology. Another application is to prove that the curvature of a line module reduces to a 2-form on the base algebra. We also extend our geometric approach to Dirac operators. We illustrate the theory on the $q$-sphere, the permutation group $S_3$ and the  bicrossproduct  model quantum spacetime with algebra $[r,t]=\lambda r$. 
 \end{abstract}

\section{Introduction}

Noncommutative differential geometry allows geometric ideas to extend to situations where the `coordinate algebra' $A$ is noncommutative\cite{Con}. The formulation of vector bundles as projective modules $E$, connections as maps $\nabla: E\to \Omega^1\tens_AE$  obeying a left Leibniz-type rule and principal bundles with quantum group fibre to which such $E$ are sections of associated bundles\cite{BrzMa:gau}, are all well known and there is a large modern literature. We follow here a constructive  approach in which the structures of differential geometry are built up layer by layer starting with the 1-forms $\Omega^1$ of which a recent review is \cite{Ma:ltcc}. As well as left connections  we will also need `bimodule connections', being left connections for which there exists a bimodule map $\sigma:E\tens_A\Omega^1\to \Omega^1\tens_AE$ with respect to which a right Leibniz rule also holds\cite{DVM,Mou}. This approach allows for the formulation of metric compatibility and is an approach particularly to noncommutative Riemannian geometry\cite{BMriem,BegMa2,BMdirac}.

In this paper we take a close look at the curvature of a connection in order to address two fundamental issues that remain open within the constructive approach. Namely we establish the two Bianchi identities and in the metric bimodule connection case the antisymmetry with respect to the metric of the Riemann tensor. From a mathematical side the covariant differentiation of morphisms will allow us to follow Chern's derivation of characteristic classes as the trace of the curvature with the result independent of the connection. This reproduces the Chern-Connes pairing between K-theory and cyclic cohomology\cite{Con} but now in a bottom up way thereby closing the gap between the Chern and the Connes approaches. Likewise from the physics side, a deeper understanding of the full Bianchi identities and the symmetries of the Riemann tensor is likely a necessary step towards the right concept of Einstein tensor in noncommutative geometries. The latter is an open problem which is plausibly key to applications of noncommutative geometry to model Planck scale or quantum gravity effects. Such effects are expected to result in noncommutative or  `quantum' spacetime of which there are now numerous models; see \cite{BegMa2,MaTao1} for two that include noncommutative Riemannian curvature.

The first of the Bianchi identities follows directly from definition of curvature, but the 2nd Bianchi identity and further symmetry or the curvature are less immediate and arise in  formulation in the context of a certain DG category ${}_A\mathcal{G}$ for noncommutative differential geometry, which in classical geometry would simply be an implicit structure behind the symmetry rules for permuting indices for covariant derivatives of tensors, but which in noncommutative geometry we need to define more explicitly in order to understand the right formulation of such symmetries. As DG categories are now central to many approaches to algebraic geometry, their natural appearance in noncommutative differential geometry is not unexpected. For an excellent reference on DG categories see \cite{dgCatKeller}. 

When Chern wrote his fundamental work on characteristic classes in differential geometry, cyclic homology and cohomology had yet to be invented.  Instead of a pairing with $K$-theory Chern's proof relied on showing that the trace of the power of the curvature was independent of the connection chosen on a particular vector bundle. With the machinery of DG categories we can follow Chern's original derivation, with one exception: Instead of cohomology classes we have to use the idea of $n$-cycles, which are also used in Connes' pairing\cite{Con}. Then the answers obtained are exactly the same as those given by Connes theory -- this can be seen by substituting in the Grassmann connection defined in terms of a matrix idempotent. It is the justification of the formula which is completely different. Of course Connes pairing with $K$-theory can make sense even when there is no sensible differential structure to calculate with, whereas the construction here is strictly differential geometric.

An outline of the paper is as follows.  We will give a formal definition of DG categories in the preliminary Section~\ref{sec2} but the idea is that in classical differential geometry the most natural maps between vector bundles with connections are the bundle maps which commute with the connections. However this is obviously too restrictive for many purposes, and is replaced by the idea of differentiating bundle maps. (A bundle map is simply a tensor with appropriate up and down indices, and is covariantly differentiated by standard rules.) This leads a notion of categories in which maps can be differentiated in an algebraic sense.

Section~\ref{sec2} specifies our geometric structures in noncommutative-algebraic terms starting with a differential graded algebra $(\Omega,\extd)$. We also define line bundles (`line modules') and the Fr\"olich map. Section \ref{bsect} constructs the DG category ${}_A\CG$ where objects are bundles with connection but morphisms include maps $E\to \Omega^n\tens_A F$ for each $n$. Theorem~\ref{the-bcat} proves the 2nd Bianchi identity and Corollary~\ref{pute} revisits the 1st Bianchi identity in our formulation. Section~\ref{secchar} contains the Chern construction with connection-independence proven in Theorem~\ref{indepconn}. Section~\ref{secbicon} switches to bimodule connections and introduces a notion of extendability that we need. Section~\ref{secAGA} introduces the bimodule version ${}_A\CG_A$ of bundles with bimodule connections, this time with tensor product over $A$ i.e.\ we show that this is a monoidal category, which we do in Secton~\ref{secAGAmon}. Corollary~\ref{yyuupp} is our Riemann curvature antisymmetry result when the connection preserves the quantum metric $g\in \Omega^1\tens_A\Omega^1$ (i.e., we prove a quantum version of the remaining classical symmetry of the Riemann tensor). 
A  second application, of ${}_A\CG_A$ is Proposition~\ref{cohommult} where we define a cup product in `sheaf cohomology' (defined as the cohomology of the noncommutative de Rham complex twisted by a flat connection). A third application, Proposition~\ref{promisedphi},  is to line bundles where we show that the curvature canonically defines a 2-form on the base of the bundle just as in classical electromagnetism.
A fourth application Corollary~\ref{curvtrace} uses a quantum metric on a module to take the trace of the curvature of a connection to give a noncommutative de Rham class in the sense of the cohomology of $(\Omega,\extd)$. We also comment on the square of the geometric Dirac operator as an extension of \cite{BMdirac}, to be developed further elsewhere.

 The paper concludes in Section~\ref{secex} with some detailed examples including the quantum cotangent bundle on the group $S_3$ of permutations of three elements, the 2-dimensional quantum spacetime model and the standard $q$-sphere from the $q$-Hopf fibration with its associated $q$-monopole bundles $E_n$ and $q$-Riemannian geometry. These are basic `test examples' of noncommutative differential geometry in our approach \cite{BMriem,BMdirac,BegMa2,Ma:spin,Ma:ltcc,Ma:bfou}. For examples like the $q$-sphere, extendability does not impose further restrictions thanks to our basic Lemma~\ref{dinnu}. Thus all the charge $n$  $q$-monopole connections and the canonical torsion-free metric compatible or  `quantum-Levi Civita' connection $\nabla_{\Omega^1}$ in \cite{Ma:spin} are extendable. The quantum antisymmetry holds and we compute the  quantum metric trace of the curvature, which classically is but now has a $q^4-1$ factor. We also exhibit our $q$-Lichnerowicz formula for the Dirac operator in \cite{Ma:spin,BMdirac}. The other two examples are at the other extreme in the sense that we show in these cases extendability of a connection on $\Omega^1$ is linked to flatness under certain homogeneity or symmetry conditions. For the bicrossproduct model spacetime we find, in particular,  a unique metric compatible extendable connection among the homogenous type studied in \cite{BegMa2}, which is flat but with torsion. For $S_3$  with its 2-cycles calculus we find exactly four metric compatible extendable connections among the ad-invariant connections studied in \cite{BMriem}. Such flat connections classically are data for the quantisation of the differential structure in Poisson-Riemannian geometry\cite{BegMa3} and it is remarkable that quantum versions arise now form our algebraic criterion of extendability.  We also exhibit the full moduli of torsion-free and cotorsion-free or `weak quantum Levi-Civita' connections and of metric compatible with-torsion connections, and illustrate one of the Bianchi identities. 
 
\section{Preliminaries}\label{sec2}
Recall that a differential graded algebra (DGA) $\Omega$ means graded components $\Omega^n$ for $n\in\N$, differentials $\extd:\Omega^n\to\Omega^{n+1}$ and an associative product $\wedge:\Omega^n\tens \Omega^m\to \Omega^{n+m}$ obeying
\[
\extd(\xi\wedge\eta)=(\extd\xi)\wedge\eta+(-1)^n\, \xi\wedge(\extd\eta)\ ,\ \extd^2=0
\]
for $\xi\in\Omega^n$ and $\eta\in\Omega^m$. We suppose that $A=\Omega^0$ is a unital algebra, that $\Omega$ is generated by $A$ under 
$\extd$ and $\wedge$, and then call $(\Omega,\extd,\wedge)$ a differential calculus for $A$. We shall work over the general field $\kk$ where possible. 

Recall that a left connection (or covariant derivative) $\nabla_E:E\to \Omega^1\tens_A E$ for a left $A$-module $E$ is a linear map obeying the left Leibniz rule,
\[
\nabla_E(a.e)=\extd a\tens e+a.\nabla_E (e)
\]
for all $a\in A$ and $e\in E$. Its curvature is the left module map 
\[
R_E=(\extd\tens\id-\id\wedge\nabla_E)\nabla_E:E\to\Omega^2 \otimes_A E.
\]
We also need right connections, which we write as $\tilde\nabla_F:F\to F\tens_A \Omega^1$ for a right $A$-module $F$, and which obey the right Leibniz rule
\[
\tilde\nabla_F(f.a)=f\tens \extd a+(\tilde\nabla_F f)a
\]
for $a\in A$ and $f\in F$. 

In the special case where we have a left connection $\nabla_{\Omega^1}$ on the bimodule $\Omega^1$ we 
have the further notion of `torsion'
\[
T=  \wedge\nabla_{\Omega^1}-\extd:\Omega^1\to \Omega^2. 
\]
We may also be interested in a right connection $\tilde\nabla_{\Omega^1}$ in which case the formula for the torsion is $\tilde T =  \wedge\tilde\nabla_{\Omega^1}+ \extd:\Omega^1\to\Omega^2$. A brief check shows that $T$ is a left module map and that $\tilde T$ is a right  module map.  A modest first result in this case, which we state for left connections (there is a similar formula for right ones) is the following.

\begin{lemma}[1st Bianchi identity]\label{1stbia} For a left connection on $\Omega^1$,
\[ \wedge R_{\Omega^1}=\extd T - (\id \wedge T)\nabla_{\Omega^1}.\]
\end{lemma}
\proof We adopt the shorthand $\nabla_{\Omega^1}\xi=\xi^1\tens\xi^2$ then
\begin{align*} \wedge R_\Omega^1(\xi)&=(\extd \xi^1)\wedge\xi^2- (\xi^1\wedge \xi^2{}^1)\wedge \xi^2{}^2\\
&=\extd(\xi^1\wedge\xi^2)+\xi^1(\extd \xi^2-\xi^2{}^1\wedge\xi^2{}^2)=\extd^2\xi+\extd T(\xi)- \xi^1\wedge T(\xi^2)
\end{align*}
for all $\xi\in \Omega^1$. We use $\extd^2=0$ and we used associativity of the wedge product. \endproof

For quantum Riemannian geometry one also has a metric $g\in \Omega^1\tens_A\Omega^1$ and in this context we define the
`cotorsion' of a left connection as
\[ {\rm co}T=(\extd\tens\id-\id\wedge\nabla_{\Omega^1})g\]
of which the vanishing is a weak notion of metric-compatibility\cite{Ma:riem,BMriem,Ma:ltcc}. For full metric compatibility we will need
more structure as we describe later in Section~\ref{secbicon}. 

Next we recall that a left $A$-module $E$ has dual the right module $E^\flat={}_A\Hom(E,A)$ (the left module maps from $E$ to $A$) with evaluation a bimodule map $\ev:E\tens E^\flat\to A$. By definition, a left finitely generated projective (left fgp for short) $A$-module $E$ has a dual basis $e^i\in E$ and $e_i\in E^\flat$ (for a finite number of indices $i$) so that 
(sum understood)
\begin{align} \label{puty}
 (\ev\tens\id) (\id\tens e_i\tens e^i)=\id:E\to E\ .
\end{align}
In the case where $E$ is also a bimodule we can define a coevaluation bimodule map $\coev:A\to E^\flat\tens_A E$ 
 so that 
$\coev(1)=e_i\tens e^i$ (sum understood).  There are formulae for the change of dual basis which read
 \begin{eqnarray} \label{changeleftbasis}
e^i=c_j(e^i).c^j,\quad e_i=c_j.e_i(c^j),\quad c^j=e_i(c^j).e^i,\quad c_j=e_i.c_j(e^i).
\end{eqnarray}
for an alternative dual basis $c^j\in E$ and $c_j\in E^\flat$. 
We shall need a result on dual connections on $E^\flat$ from \cite{BMriem}: 
For the left fgp module $E$ with left connection $\nabla_E$, there is a unique right connection $\tilde\nabla_{E^\flat} : E^\flat \to  E^\flat \otimes_{A} \Omega^{1} $ 
given by
\begin{eqnarray*}
\tilde\nabla_{E^\flat}(\alpha)=e_j\tens \extd\,\ev(e^j\tens \alpha) - 
e_j\tens (\id\tens\ev_E)(\nabla_E(e^j)\tens \alpha).      
\end{eqnarray*}
so that
\begin{eqnarray*}
\extd\circ\ev=(\id\tens\ev)(\nabla_E\tens\id)+(\ev\tens\id)(\id\tens\tilde\nabla_{E^\flat}):E \tens E^\flat\to\Omega^1 .
\end{eqnarray*}
Now we apply the covariant derivatives to the dual basis to get
\begin{align*}
&( \tilde\nabla_{E^\flat}\tens\id+\id\tens\nabla_E)(e_i\tens e^i)   \cr
&\ =\ e_j\tens \extd\,\ev(e^j\tens e_i) \tens e^i   - 
e_j\tens (\id\tens\ev_E)(\nabla_E(e^j)\tens e_i)\tens e^i +e_i \tens\nabla_E(e^i) \cr
&\ =\ e_j\tens \extd\,\ev(e^j\tens e_i) \tens e^i
\end{align*}
where we have used (\ref{puty}). Now setting $P_{ji}=\ev(e^j\tens e_i)$ we have
\[
( \tilde\nabla_{E^\flat}\tens\id+\id\tens\nabla_E)(e_i\tens e^i)=e_j\tens\extd P_{ji} \tens e^i=e_q\tens P_{qj}.\extd P_{ji}.P_{ik} \tens e^k
\]
and applying $\extd$ to the matrix product $P^2=P$ we get $P.\extd P.P=0$, so 
\begin{align} \label{oscat}
( \tilde\nabla_{E^\flat}\tens\id+\id\tens\nabla_E)(e_i\tens e^i)=0.
\end{align}
We also recall the definition of a line module, which is the direct noncommutative generalisation of a line bundle, over a unital algebra $A$. This mean a left fgp $A$-bimodule $L$ for which the bimodule maps of coevaluation $\coev:A\to L^\flat\tens_A L$
and evaluation $\ev:L\tens_A L^\flat\to A$ are isomorphisms. Such modules were originally introduced as Morita contexts in algebraic $K$-theory (see \cite{BassK}), and later as geometric objects in \cite{NCline}.  A standard result we will need later (from e.g.\ \cite{BassK}) is:

\begin{lemma} \label{inventelement}
Let $L$ be a line $A$-module. If   $T:L\to L\tens_A F$ is a left $A$-module map for some left $A$-module $F$ then
there is a unique $f\in F$ such that $T(e)=e\tens f$.  If  $S:L\to E\tens_A L$ is a right $A$-module map for some right $A$-module $E$ then
there is a unique $g\in E$ so that $S(e)=g\tens e$. 
\end{lemma}

It is shown in \cite{fropicard} that associated to a line module $L$ is a unique  unital algebra automorphism $\Phi_L$ of the centre $Z(A)$ (we will call it the Fr\"ohlich map) such that  
\[ z.e=e.\Phi_L(z)\]
for all $e\in L$ and  $z\in Z(A)$. This map
depends only on the isomorphism class of the bimodule and obeys  $\Phi_{L\tens M}=\Phi_M\circ\Phi_L$. 

We will also need a little category theory. For any field $\kk$, we define (see \cite{catmitchell}) a $\kk$-category\index{$\kk$-category} to be an additive category with a ring homomorphism from $\kk$ to the natural transformations from the identity functor to itself. An example of a $\kk$-category is the category of vector spaces and linear maps over $\kk$, where $\kk$ acts by multiplication on each object. Now we define a DG category (see \cite{dgCatKeller}):

\begin{definition} \label{dgcatdef}
A differential graded category\index{differential graded category} (or DG category) is a $\kk$-category where each $\Mor(X,Y)$ is a cochain complex (in a weak sense where we do not require the differential squared to be zero) of $\kk$-vector spaces, and composition $\circ:\Mor(Y,Z)\tens \Mor(X,Y)\to \Mor(X,Z)$ is a map of cochain complexes.
\end{definition}

In what follows we use an order of products which is natural for modules and connections, even if it is not standard for DG categories. 

\section{The DG category ${}_A\mathcal{G}$} \label{bsect}
Given two left $A$-modules with connection $(E,\nabla_E)$ and $(F,\nabla_F)$ a left module map $\phi:E\to F$ is said to intertwine the connections if 
\[
\nabla_F\circ\phi=(\id\tens\phi)\circ\nabla_E:E\to\Omega^1\otimes_A F.
\]
If we generalise the maps which we consider to all left module maps then it is useful to define the derivative of $\phi:E\to F$ to be
\[
\doublenabla(\phi)=\nabla_F\circ\phi-(\id\tens\phi)\circ\nabla_E:E\to\Omega^1\otimes_A F
\]
so $\doublenabla(\phi)=0$ precisely when $\phi$ intertwines the connections. Note that $\doublenabla(\phi)$ is a left module map when $\phi$ is. It is convenient to form a category with objects left modules with connections, and to grade the morphisms from $E$ to $F$ by $\N$ so that $\phi:E\to F$ is in $\Mor_0(E,F)$ and $\doublenabla(\phi)$ is in $\Mor_1(E,F)$. Our basic example of a DG category for noncommutative geometry will be based on the following category with graded morphisms.

\begin{definition} \label{kjhxdhvc}
The category ${}_A\mathcal{G}$\index{category!${}_A\mathcal{G}$} is defined with objects
$(E,\nabla_E)$  left $A$-modules $E$ 
with left connections. Morphisms are graded by $\mathbb{N}$ (including
$0\in\mathbb{N}$) with $\psi\in \mathrm{Mor}_n((E,\nabla_E),(F,\nabla_F))$ a 
 left module map $\psi:E\to \Omega^n\tens_A F$.
  Composition of morphisms is given by the formula
\begin{eqnarray*}
\phi\circ\psi\,=\,(\id\wedge\phi)\psi:E\to \Omega^{n+m}\otimes_A G 
\end{eqnarray*}
where $\phi\in 
\mathrm{Mor}_m((F,\nabla_F),(G,\nabla_G))$.
\end{definition}

 To give the formula for derivatives of morphisms in all grades, it will be convenient to extend the 
left connection $\nabla_E : E\to\Omega^1 \tens_A E$ to
\begin{eqnarray} \label{extendconnection}
\nabla_E{}^{[n]}:\Omega^n 
\otimes_A E\to \Omega^{n+1} \otimes_A E\, , \quad \xi\tens e
\mapsto 
\extd\xi\tens e+(-1)^n\xi\wedge \nabla_E e 
\end{eqnarray}
for $n\ge 0$, the $n=0$ case giving $\nabla_E$. Using this we can write the curvature $R_E$ of the connection $(E,\nabla_E)$ as
\[
R_E=\nabla_E{}^{[1]}\,\nabla_E{}^{[0]}:E\to\Omega^2\otimes_A E.
\]

\begin{theorem} \label{the-bcat}
For $\psi:E\to \Omega^n\tens_A F$  an $n$-morphism in ${}_A\mathcal{G}$ we define its differential
 \[ 
\doublenabla(\psi)\,=\,\nabla_F{}^{[n]}\circ \psi-(\id\wedge\psi)\nabla_E:
E\to \Omega^{n+1}\otimes_A F\]
as an $n+1$-morphism. 
The category ${}_A\mathcal{G}$ with this differential is a DG category. 
 If  $\phi:F\to \Omega^m\tens_A G$ is an $m$-morphism then 
\begin{eqnarray*}
\doublenabla(\phi\circ\psi) \,=\, \phi\circ\doublenabla(\psi)+(-1)^n\doublenabla(\phi)\circ\psi.
\end{eqnarray*}
Moreover,

{\rm (1)}  (2nd Bianchi) The curvature $R_E\in \mathrm{Mor}_2((E,\nabla_E),(E,\nabla_E))$ has $\doublenabla(R_E)=0$;

{\rm (2)}  For all morphisms $\psi:E\to F$, 
$\doublenabla(\doublenabla(\psi))=R_F\circ \psi-\psi\circ R_E$.
\end{theorem}
\proof We first show that $\doublenabla(\psi)$ is a left module map. For $a\in A$, $e\in E$,
and  $\psi(e)=\xi\tens f$ (summation implicit), 
\begin{align*}
&\kern-10pt\nabla^{[n]}_F\circ \psi(a.e)-(\id\wedge\psi)\nabla_E(a.e) = 
\nabla^{[n]}_F(a.\xi\tens f)-(\id\wedge\psi)(\extd a\tens e+a.\nabla_E(e)) \cr
=&\  \extd a\wedge \xi\tens f+
 a.\nabla^{[n]}_F(\xi\tens f)-\extd a\wedge\psi( e)+a.\nabla_E(e) 
 = a.(\nabla^{[n]}_F\circ \psi(e)-(\id\wedge\psi)\nabla_E(e)).
\end{align*}
Next setting $\phi(f)=\kappa\tens g$ (summation implicit),
\begin{align*}
&\kern-20pt\nabla_G{}^{[n+m]}(\id\wedge\phi)\psi(e) =  \nabla_G{}^{[n+m]}(\xi\wedge\phi(f)) 
 = \nabla_G{}^{[n+m]}(\xi\wedge\kappa\tens g)
 \cr =&\ \extd(\xi\wedge\kappa)\tens g+(-1)^{n+m}\xi\wedge\kappa\wedge\nabla_G g 
  = \extd\xi\wedge\phi(f)+(-1)^n\xi\wedge\nabla_G{}^{[m]}(\phi(f))\cr
  =&\   (\id\wedge\phi)\nabla_F{}^{[n]}\psi(e) + (-1)^n\xi\wedge\doublenabla(\phi)f\cr
  =&\  (\id\wedge\phi)\doublenabla(\psi)e
 +(\id\wedge(\id\wedge\phi)\psi)\nabla_E e
  + (-1)^n\xi\wedge\doublenabla(\phi)f  .
\end{align*}
For (1), 
\begin{eqnarray*}
\doublenabla(R_E) &=& \nabla_E{}^{[2]}\circ R_E - (\id\wedge R_E)\nabla_E \cr
&=&  \nabla_E{}^{[2]}\circ  \nabla_E{}^{[1]}\circ \nabla_E - (\id\wedge (\nabla_E{}^{[1]}\circ \nabla_E))\nabla_E .
\end{eqnarray*}
Now set $\nabla_E e=\xi\tens h$ and $\nabla_E h=\eta\tens f$ (summation implicit),
\begin{align*}
&\kern-20pt \doublenabla(R_E)(e) = 
  \nabla_E{}^{[2]}( \nabla_E{}^{[1]}(\xi\tens h)) - \xi\wedge (\nabla_E{}^{[1]}\circ \nabla_E)h \cr
  =&\  \nabla_E{}^{[2]}(\extd\xi\tens h - \xi\wedge\nabla_E h) - \xi\wedge (\nabla_E{}^{[1]}\circ \nabla_E)h \cr
=&\ -\,\extd(\xi\wedge\eta)\tens f+\extd\xi\tens \nabla_E h - \xi\wedge\eta\wedge\nabla_E f  -  \xi\wedge \extd\eta\tens f + \xi\wedge \eta\wedge\nabla f = 0.
\end{align*}
To prove (2) we begin by showing that for all $n\ge 0$, 
\[\nabla_E{}^{[n+1]}\circ \nabla_E{}^{[n]}=\id\wedge R_E:\Omega^n 
\otimes_A E\to \Omega^{n+2} \otimes_A E.
\]
Putting $\nabla_E e=\xi \tens h $ again for convenience, and for all $\omega\in \Omega^n$, 
\begin{align*}
&\kern-20pt\nabla_E{}^{[n+1]}( \nabla_E{}^{[n]}(\omega\tens e)) =
\nabla_E{}^{[n+1]}(\extd\omega\tens e+(-1)^n\omega\wedge \nabla_E e)\cr
=&\
\nabla_E{}^{[n+1]}(\extd\omega\tens e+(-1)^n\omega\wedge \xi \tens h ) \cr
=&\ (-1)^{n+1}\extd\omega\wedge \nabla_E e+
(-1)^n\extd\omega\wedge \xi \tens e 
+\omega\wedge \extd\xi \tens h  -  \omega\wedge \xi \wedge \nabla_E h  \cr
=&\ \omega\wedge(\extd\xi \tens h 
- \xi \wedge \nabla_E h )\,=\,\omega\wedge R_E(e).
\end{align*}
To complete the proof of (2),
\begin{align*}
\doublenabla(\doublenabla(\psi))=&\ \nabla_F{}^{[n+1]}\circ \doublenabla(\psi)-(\id\wedge\doublenabla(\psi))\nabla_E \cr
=&\  \nabla_F{}^{[n+1]}\circ \nabla_F{}^{[n]}\circ \psi-\nabla_F{}^{[n+1]}\circ (\id\wedge\psi)\nabla_E
-(\id\wedge(\nabla_F{}^{[n]}\circ \psi))\nabla_E \cr
&\ +\ (\id\wedge\psi)(\id\wedge\nabla_E)\nabla_E \cr
=&\  (\id\wedge R_F) \psi-\nabla_F{}^{[n+1]}\circ (\id\wedge\psi)\nabla_E
-(\id\wedge(\nabla_F{}^{[n]}\circ \psi))\nabla_E \cr
&\ -\ (\id\wedge\psi)R_E
+\ (\id\wedge\psi)(\extd\wedge\id)\nabla_E .
\end{align*}
and if $\nabla_E e=\xi\tens h$ and $\psi(h)=\eta\tens f$ then
\begin{eqnarray*}
\nabla_F{}^{[n+1]}\circ (\id\wedge\psi)\nabla_E(e) &=& \nabla_F{}^{[n+1]}(\xi\wedge\eta\tens f)\cr
&=& \extd \xi\wedge\eta\tens f - \xi\wedge\extd\eta\tens f-(-1)^n \xi\wedge\eta\wedge\nabla_F f,\cr
(\id\wedge(\nabla_F{}^{[n]}\circ \psi))\nabla_E(e) &=& \xi\wedge\nabla_F{}^{[n]}(\eta\tens f) \cr
&=& \xi\wedge\extd\eta\tens f+(-1)^n \xi\wedge\eta\wedge\nabla_F f. \qquad\square \end{eqnarray*}

In this result (1) becomes the usual 2nd Bianchi identity in the classical case when this is viewed in $\Omega^3$, while (2) tells that the curvatures form an obstruction to the cochain complex on each morphism space having differential $\doublenabla$ that squares to zero. This means that each $(\mathrm{Mor}(E,E),\doublenabla,R_E)$ in ${}_A\mathcal{G}$ forms a curved DGA as described in \cite{blockcat}.  We summarise the definition of the category ${}_A\mathcal{G}$ in a table:

\medskip
\noindent
{\renewcommand{\arraystretch}{1}
\begin{tabular}{c|c|c|c}Name & Objects &   $\Mor_n((E,\nabla_E),(F,\nabla_F))$ & derivative  \\
\hline
${}_A\mathcal{G}$\index{category!${}_A\mathcal{G}$}  & $(E,\nabla_E)$ &   $\phi:E\to\Omega^n \otimes_A F$  & 
$\doublenabla:\Mor_n(E,F)\to \Mor_{n+1}(E,F)$   \\
 & left modules \&   &   left module maps & 
 $ \doublenabla(\phi)\,=\,\nabla_F{}^{[n]}\circ \phi-(\id\wedge\phi)\nabla_E$  \\
  &  \ left connections  &  
\end{tabular}}

\medskip Now we show how some almost tautological morphisms in ${}_A\mathcal{G}$ can have a nice geometric interpretation, including a version of the first Bianchi identity\index{Bianchi identity!first}.

\begin{example} \label{firbia} Let $A$ be an algebra with differential structure and consider $(A,\extd)$ as a bimodule with connection. Suppose that $\Omega^1$ has a left connection $\nabla$ and consider the morphism
$\tau\in \mathrm{Mor}_1(\Omega^1 ,A)$ given by $\xi\mapsto \xi\tens_A 1$. Then $(\doublenabla \tau)(\xi) = \extd\xi\tens 1 -(\id\wedge\tau)\nabla\xi\,=\,(\extd\xi-\wedge\nabla\xi)\tens 1$ showing that 
$\doublenabla \tau=-T_\nabla\in  \mathrm{Mor}_2(\Omega^1 ,A)$ where $T_\nabla$ is the
 torsion.
 \end{example} 
 
 There is a corollary of Theorem~\ref{the-bcat} which will give us a more conceptual way of thinking about the 1st Bianchi identity. We note that classically the vector fields are dual to the 1-forms, and this continues to provide a good definition of vector field in noncommutative geometry (as opposed to the idea of identifying vector fields with derivations). Of course in noncommutative geometry we have to choose a side. 
 
 \begin{corollary}[1st Bianchi identity revisited] \label{pute}
 Suppose that $\Omega^1$ is right finitely generated projective with dual $\cX^R$ (the right vector fields) and suppose $\cX^R$ has a left connection with $(\cX^R,\nabla_{\cX})\in {}_A\mathcal{G}$. Then for the coevaluation map
$\coev:A\to\Omega^1\tens_A \cX^R$,
\[
\doublenabla(\doublenabla(\coev))=R_\cX\circ \coev\in \mathrm{Mor}_3(A,\cX^R).
\]
\end{corollary} 
\proof Under our assumptions,
$\coev:A\to\Omega^1\tens_A \cX^R$ is in  $\mathrm{Mor}_1(A,\cX^R)$ and from Theorem~\ref{the-bcat}
we have the 1st Bianchi identity as stated.   \endproof

At first sight this looks nothing like the first Bianchi identity but remember that
if we set $E$ to be the left fgp module $\cX^R$ then its dual $E^\flat$ is $\Omega^1$ so we have a corresponding right connection $\tilde\nabla=\tilde\nabla_{\Omega^1}$ on $\Omega^1$ such that (\ref{oscat}) becomes
\begin{align}\label{osdog}
(\tilde\nabla_{\Omega^1}\tens\id + \id\tens\nabla_{\cX})\coev(1)=0 \in\Omega^1\otimes_A \Omega^1\otimes_A \cX^R.
\end{align}

\begin{corollary}\label{pute2} If the torsion $\tilde T$ of $\tilde\nabla$ is a bimodule map then the 1st Bianchi identity in Corollary~\ref{pute}  becomes 
\[
\doublenabla((\tilde T\tens \id)\coev)=(\id\wedge R_\cX) \coev:A\to \Omega^3 \otimes_A \cX^R.
\]
\end{corollary}
\proof Applying $\wedge\tens\id$ to (\ref{osdog}) gives
\begin{align*}
((\wedge \tilde\nabla_{\Omega^1})\tens\id + \id\wedge\nabla_{\cX}) \coev(1)=0 \in  \Omega^2\otimes_A \cX^R
\end{align*}
which we rewrite as
\[
(\extd\tens\id +(\wedge\tilde\nabla_{\Omega^1})\tens\id) \coev(1)=
(\extd\tens\id - \id\wedge\nabla_{\cX})\coev(1). 
\]
Hence if $\tilde T=\extd+(\wedge\tilde\nabla_{\Omega^1}):\Omega^1\to\Omega^2$, i.e.\ the (right handed) torsion of $\tilde\nabla$, is a left (and hence bi-)module map then 
\[
\doublenabla(\coev) = (\tilde T\tens \id)\coev \in \mathrm{Mor}_2(A,\cX^R)
\]
and we then use Corollary~\ref{pute} to obtain the stated result.
\endproof 

The form version of this is the right handed version of Lemma~\ref{1stbia}, for example $\wedge\tilde R=0$ in the case of zero torsion, but the above gives the Bianchi identity more conventionally expressed in terms of the curvature on vector fields. In the classical case with zero torsion and using the usual conventions for the Riemann tensor, our result in Corollary~\ref{pute2} reads 
\[
R^a{}_{bcd}\,\extd x^b \wedge \extd x^c \wedge \extd x^d \tens\tfrac{\pd}{\pd x^a}=0
\]
which is the familiar classical first Bianchi identity in the absence of torsion.

\section{Characteristic classes}\label{secchar}
 Chern defined characteristic classes of a vector bundle in terms of the trace of powers of the curvature \cite{Chernclass}. To follow this construction literally in noncommutative geometry would require a formulation of trace, which is problematic without additional structures such as braiding. However, we can use a rather more limited version of trace. 

\begin{lemma} \label{limitedtrace}
Let $F$ be an $A$-bimodule, $E$ a left finitely generated projective module and a linear map $\phi:F\to \kk$ have the trace property $\phi(a.f)=\phi(f.a)$ for all $a\in A$, $f\in F$. Then a left module map $\theta:E\to F\tens_A E$ has a well defined trace
\begin{eqnarray*}
\trace_\phi(\theta)=\phi(  (\id\tens\ev)  (\theta e^i\tens e_i)) \in \kk  .
\end{eqnarray*}
independently of the choice of dual basis $e^i$ of $E$ and $e_i$ of $E^\flat$.
\end{lemma}
\proof We need to show that the formula does not depend on the choice of dual basis.
If we take a different dual basis  $c_j,c^j$ of $E^\flat,E$ then
\begin{align*}
&\kern-10pt\phi(  (\id\tens\ev)  (\theta(c^j)\tens c_j)) =
\phi(  (\id\tens\ev)  (\theta(e_i(c^j).e^i)\tens e_k.c_j(e^k))) \cr
=&\ 
\phi(  e_i(c^j).(\id\tens\ev)  (\theta(e^i)\tens e_k.c_j(e^k))) 
= 
\phi(  (\id\tens\ev)  (\theta(e^i)\tens e_k.c_j(e^k)\,e_i(c^j))) \cr
=&\ 
\phi(  (\id\tens\ev)  (\theta(e^i)\tens e_k.e_i(c_j(e^k)\,c^j))) 
=
\phi(  (\id\tens\ev)  (\theta(e^i)\tens e_k.e_i(    e^k      ))) \cr
=&\ 
\phi(  (\id\tens\ev)  (\theta(e^i)\tens e_i   ))
\end{align*}
on using the change of basis formulae in (\ref{changeleftbasis}) and summing over $i,j,k$. \eproof

The linear maps we wish to use in this result are the $n$-cycles as in \cite{Con}. Given a differential calculus $\Omega$ on $A$, 
an $n$-cycle is a linear map $\int:\Omega^n\to \kk$
satisying the conditions
\begin{eqnarray} \label{ncycledef}
\int \extd\xi=0,\quad 
 \int \omega\wedge \rho=(-1)^{|\omega||\rho|}\int\rho\wedge\omega
\end{eqnarray}
for all $\xi\in\Omega^{n-1}$ and all forms $\omega,\rho$ with sum of degrees $|\omega|+|\rho|=n$.
We now combine $n$-cycles with the derivatives of morphisms in ${}_A\mathcal{G}$.

\begin{proposition} \label{limitedtraced} 
Let $E$ be an $A$-bimodule which is finitely generated projective as a left module, let $\nabla_E$ be a left connection and $\int$ an $n+1$-cycle on $\Omega$. If  $(E,\nabla_E)\in {}_A\mathcal{G}$ then any $\theta\in \mathrm{Mor}_n(E,E)$ has 
$\trace_{\int}(\doublenabla(\theta))=0$.
\end{proposition}
\proof   Take a dual basis $e^i\in E$ and $e_i\in E^\flat$, and another copy
$e^j\in E$ and $e_j\in E^\flat$. Using the dual connection, we have
\begin{align*}
&\kern-20pt\trace_{\int}(\doublenabla(\theta))=\ {\int}(\id\tens\ev_E)  (\doublenabla(\theta)(e^i)\tens e_i ) \cr
=& {\int} (\id\tens\ev_E)\big((\extd\tens\id+(-1)^n(\id\wedge\nabla_E))\theta(e^i)\tens e_i- (\id\wedge\theta)\nabla_E(e^i)\tens e_i\big) \cr
=&\ {\int} \extd (\id\tens\ev_E)  (\theta(e^i)\tens e_i) -(-1)^n
{\int}(\id\wedge\ev_E\wedge\id)  (\theta(e^i)\tens \tilde\nabla_{E^\flat}(e_i))\cr
&\  -\, {\int}  (\id\tens\ev_E) ( (\id\wedge\theta)\nabla_E(e^i)\tens e_i) .
\end{align*}
Using that ${\int}$ is an $n+1$ cycle and the shorthand $\tilde\nabla_{E^\flat}(e_i)=e_j\tens \eta_{ji}$ we have\begin{eqnarray*}
&& (-1)^n
{\int}(\id\wedge\ev_E\wedge\id)  (\theta(e^i)\tens \tilde\nabla_{E^\flat}(e_i))   
= {\int} (\id\tens\ev_E)(\eta_{ji}\wedge \theta(e^i)\tens e_j)
\end{eqnarray*}
while using the explicit formula for $\eta_{ji}$ in the dual connection and the shorthand $\nabla_E(e^j)=\kappa\tens f$  (summation implicit),
\begin{align*}
\eta_{ji}\wedge \theta(e^i) =&\ \extd\,\ev(e^j\tens e_i) \wedge \theta(e^i) 
- (\id\tens\ev_E)(\nabla_E(e^j)\tens e_i)\wedge \theta(e^i) \cr
=&\ \extd\,\ev(e^j\tens e_i) \wedge \theta(e^i) 
- \kappa\ \ev_E(f\tens e_i)\wedge \theta(e^i) \cr
=&\ \extd\,\ev(e^j\tens e_i) \wedge \theta(e^i) 
- \kappa\, \wedge \theta(\ev_E(f\tens e_i)\,e^i) \cr
=&\ \extd\,\ev(e^j\tens e_i) \wedge \theta(e^i) 
- \kappa\, \wedge \theta(f) \cr
=&\ \extd\,\ev(e^j\tens e_i) \wedge \theta(e^i) 
- (\id\wedge \theta)\nabla_E(e^j) .
\end{align*}
Substituting our these expressions and setting $P_{ji}=\ev(e^j\tens e_i)$ so that $e^i=P_{ik}e^k$ and $e_j=e_m\,P_{mj}$, we arrive at
\begin{align*}
\trace_{\int}(\doublenabla(\theta))
=&\ - \, {\int} (\id\tens\ev_E)(\extd\,\ev(e^j\tens e_i) \wedge \theta(e^i) \tens e_j)\cr
=&\ - \, {\int}(\id\tens\ev_E)(P_{mj}.\extd P_{ji}.P_{ik} \wedge \theta(e^k) \tens e_m)  .
\end{align*}
But applying $\extd$ to $P^2=P$ tells us that $P.\extd P.P=0$.\eproof

To apply these results to characteristic classes, we use powers of the curvature, and for $(E,\nabla_E)\in {}_A\mathcal{G}$ it is natural to take these as 
 $R_E{}^2=R_E{}\circ R_E{}$, $R_E{}^3=R_E{}\circ R_E{} \circ R_E{}$ etc. Now from Theorem~\ref{the-bcat} we have
$\doublenabla(R_E{}^n)=0$. 
The next result is a noncommutative version of a standard classical result at the heart of Chern's construction; that the cohomology class of the trace of the curvature is independent of the connection:

\begin{theorem}  \label{indepconn}
Suppose  that $E$ is a finitely generated projective left $A$-module, and that $\int$ is a $2n$-cycle on the differential calculus $\Omega$.
Then for any left covariant derivative $\nabla_E$ and any $n$, the value of $\trace_{\rm \int}(R_E{}^n)$ is independent of the choice of $\nabla_E$. 
\end{theorem}
\proof The set of left connections on $E$ is an affine space, so any two left connections can be connected by a straight line path. Take such a path, with connection $\nabla_E{}^t$ parameterised by $t\in\R$, and define
$\dot\nabla_E{}:E\to \Omega^1 A\tens_A E$ by
\begin{eqnarray*}
\dot\nabla_E{}(e)\,=\,\frac{\extd \nabla_E{}^t(e)}{\extd t}.
\end{eqnarray*}
(As we are moving along a straight line interpolating two connections, $\dot\nabla_E$ does not depend on $t$.)
For all $t$, we have
$
\nabla_E{}^t(a.e)\,=\,\extd a\tens e +a.\nabla_E{}^t(e)
$ from which it follows that $\dot\nabla_E{}$ is a left module map. Differentiating
$R_E{}^t\,=\,(\nabla_E{}^t)^{[1]}\circ\nabla_E{}^t$ (which is a quadratic function of $t$), we obtain
\begin{eqnarray*}
 \dot R_E{}^t\,=\, -(\id\wedge\dot\nabla_E{})\nabla_E{}^t+(\nabla_E{}^t)^{[1]}\dot\nabla_E{}\,=\,\doublenabla^t(\dot\nabla_E{}).
\end{eqnarray*}
Now omitting the explicit $t$ dependence for clarity, by Theorem~\ref{the-bcat},
\begin{align*}
\tfrac{\extd}{\extd t} R_E{}^n=&\ \dot R_E\circ R_E{}^{n-1} +R_E\circ \dot R_E\circ R_E{}^{n-2}+\dots
+R_E{}^{n-1}\circ \dot R_E  \cr
=&\ \doublenabla(\dot\nabla_E)\circ R_E{}^{n-1} +R_E\circ \doublenabla(\dot\nabla_E)\circ R_E{}^{n-2}+\dots
+R_E{}^{n-1}\circ \doublenabla(\dot\nabla_E) \cr
=&\ \doublenabla(\dot\nabla_E\circ R_E{}^{n-1} +R_E\circ \dot\nabla_E\circ R_E{}^{n-2}+\dots
+R_E{}^{n-1}\circ \dot\nabla_E).
\end{align*}
Then 
$\tfrac{\extd}{\extd t}\mathrm{Tr}{}_{\int}(R_E{}^n)=0$ by Proposition~\ref{limitedtraced}. 
\eproof

We still need to know that  left connections actually exist for a reasonable class of modules. 

\begin{example} \label{grassconn}
Given $E$ an fgp left $A$-module and a choice of dual bases $e^i\in E$ and $e_i\in E^\flat$, 
 we construct a left covariant derivative $\nabla_E:E\to\Omega^1\tens_A E$, called the Grassmann connection\index{Grassmann connection},   by
\begin{eqnarray*}
\nabla_E e = \sum_i\extd(e_i(e))\tens e^i,\quad e\in E.
\end{eqnarray*}
If $e_j,e^j$ is an independent dual basis pair, the curvature can be written
\begin{eqnarray*}
R_E e = -\sum_{i,j}\extd(e_i(e))  \wedge \extd(e_j(e^i))    \tens e^j.
\end{eqnarray*}
Writing $P_{ij}=\ev(e^i\tens e_j)=e_j(e^i)$ and using $P_{ij}\,e^j=e^i$, we have
\begin{eqnarray} 
\nabla_E e^i = \sum_i\extd P_{ij}.P_{jk}\tens e^k,\quad
R_E e^i = -\sum_{j,k,m} \extd P_{ij}\wedge \extd P_{jk}.P_{km}\tens e^m. \label{grassformula}
\end{eqnarray}
\end{example}

In fact we have a connection on any projective module, but we would not have such a simple formula in terms of projection matrices in the general case. (By the Cuntz-Quillen theorem
\cite{CQ} a module is projective if and only it has a connection for the \textit{universal} calculus.) Now, we can take the expression for the curvature in Example~\ref{grassconn}, insert it into the formula for the trace involving the cycle in 
Theorem~\ref{indepconn}, and observe that we have just the same formula for invariants associated to projection matrices as that given by Connes \cite{Con}. To stress the point again, our derivation is completely different -- Connes' construction was via cyclic theory, and we followed Chern's construction. We shall return to the meaning of these invariants later in Section~\ref{secAGA}.

\section{The DG category ${}_A\CG_A$}\label{secAGA}

To proceed further we will need our connections to be bimodule connections as these allow tensor products, and therefore a decent notion
of metric compatibility with respect to a metric $g\in \Omega^1\tens_A\Omega^1$. We start with a technical
extendability property, then define the DG category ${}_A\CG_A$, and then prove it is monoidal under tensor products
of bimodules over $A$. We assume throughout that $A$ is equipped with a differential calculus $(\Omega,\extd,\wedge)$.

\subsection{Bimodule connections and extendability} \label{secbicon}

For bimodules over a noncommutative algebra it turns out that asking for the left and right Leibniz rules to hold for a connection simultaneously is too strong a condition, and instead we use the following definition going back to \cite{DVM,Mou}. As mentioned in the introduction, a left {\em bimodule connection} $(\nabla_E,\sigma_E)$ on a bimodule $E$ is a
left connection $\nabla_{E} : E \to \Omega^{1} \otimes_{A} E$  such that there exists a bimodule map
$\sigma_E:E\tens_A\Omega^1 \to \Omega^1 \tens_A E$ obeying, for $e\in E$ and $a\in A$,
\begin{eqnarray*}
\nabla_E(e.a)=\sigma_E(e\tens\extd a)+(\nabla_E e).a.
\end{eqnarray*}
Note that if it exists $\sigma$ is uniquely determined so this is a property that some left connection have rather than additional data.
The most important feature of bimodule connections is that we can tensor product them. Thus, if  $(E,\nabla_E,\sigma_E)$ is a (left) bimodule connection on the bimodule $E$ and that
 $(F,\nabla_F)$ is a left connection on the left module $F$. Then we can define a tensor product connection on $E\tens_A F$ by \cite{BDMS}.
 \begin{align}\label{tens45}
\nabla_{E\tens F} =&\ \nabla_E\tens
\id+(\sigma_E\tens\id)(\id\tens\nabla_F): E\otimes_A F\to \Omega^1 \otimes_A E\otimes_A F.
\end{align}
 Further if $F$ is a bimodule with $(F,\nabla_F,\sigma_F)$ a (left) bimodule connection then we have a tensor product bimodule connection
$(E\tens_A F,\nabla_{E\tens F},\sigma_{E\tens F})$, where
\begin{align*}
\sigma_{E\tens F} =&\ (\sigma_{E}\tens
\id)(\id\tens \sigma_{F}):E\otimes_A F\otimes_A  \Omega^1  \to \Omega^1 \otimes_A E\otimes_A F.
\end{align*}
This is key to applications in noncommutative Riemannian geometry\cite{BMriem,BegMa2} as it means that a bimodule connection on $\Omega^1$ acts on the metric tensor $g\in \Omega^1\tens_A\Omega^1$ so that we can write $\nabla g=0$ (otherwise for general left connections one has only the weaker notion of cotorsion free). It also means that we have the monoidal categories  ${}_A\mathcal{E}_A $, ${}_A\mathcal{EI}_A$ as in \cite{BMriem} where the former has objects which are bimodules with bimodule connection and the latter is the same with $\sigma$ invertible. Morphisms were defined as bimodule maps which intertwine the connections.   

Before generalising ${}_A\CE_A$, we need to specify an additional possible property of bimodule covariant derivatives, which is the ability to swap the order of $n$-forms and modules in the same way that the definition of a bimodule connection incorporates swapping $1$-forms and modules.

\begin{definition} \label{extendableconn} 
An $A$-bimodule with left bimodule connection $(E,\nabla_E,\sigma_E)$ is called {\em{extendable}}\index{extendable bimodule connection} if
 $\sigma_E:E\tens_A\Omega^1 \to \Omega^1  \tens_A E$
extends for all $n\ge 1$ to $\sigma_E:E\tens_A\Omega^n \to \Omega^n  \tens_A E$ 
such that \begin{eqnarray*}
&&(\wedge\tens\id)(\id\tens\sigma_E)(\sigma_E\tens\id) = \sigma_E(\id\tens\wedge):E\otimes_A\Omega^n 
\otimes_A\Omega^m  \to \Omega^{n+m}  \otimes_A E.
\end{eqnarray*}
for all $m\ge 1$. We include the case $m,n=0$ with $\sigma_E(e\tens 1)=1\tens e$. 
\end{definition}

As $\Omega$ is generated by $A,\extd A$ the extended $\sigma_E$ (should they exist) are uniquely determined by $\sigma_E$ on the 1-forms. 
There is a useful property which holds when both extendability and curvature being a right module map hold:

\begin{lemma} \label{mor2help}
Suppose that $(E,\nabla_E,\sigma_E)$ is an extendable left $A$-bimodule connection. Then the curvature $R_E$ is a right module map if and only if 
\begin{eqnarray*}
&&\nabla_E{}^{[n]}\sigma_E = (\id\wedge\sigma_E)(\nabla_E\tens\id)+\sigma_E(\id\tens\extd) : E\otimes_A\Omega^n
\to \Omega^{n+1}\otimes_A E.
\end{eqnarray*}
holds for $n=1$. In this case it also holds for all $n\ge 0$ and we also have
 \begin{eqnarray*}
(\id\wedge R_E)\sigma_E\,=\, (\id\wedge\sigma_E)(R_E\tens\id)
   : E\otimes_A\Omega^n \to \Omega^{n+2}  \otimes_A E\ .
\end{eqnarray*}
\end{lemma}
\proof  For $a\in A$ and $e\in E$ and writing $\nabla e=\eta\tens f$ (summation implicit), we have $\nabla_E(e.a) =\ \sigma_E(e \tens \extd a) + \eta\tens f.a$ and 
\begin{align*}
R_E(e.a) =&\ (\extd\tens\id-\id\wedge\nabla_E)\sigma_E(e \tens \extd a) 
+ \extd\eta\tens f.a -  \eta\wedge\nabla_E( f.a)  \cr
=&\ (\extd\tens\id-\id\wedge\nabla_E)\sigma_E(e \tens \extd a)  
+ \extd\eta\tens f.a -  \eta\wedge(\nabla_E f).a   -  \eta\wedge\sigma_E(f\tens\extd a)\cr
=&\ (\extd\tens\id-\id\wedge\nabla_E)\sigma_E(e \tens \extd a)  
+ R_E(e).a    -  \eta\wedge\sigma_E(f\tens\extd a)\ .
\end{align*}
Thus $R_E$ is a right module map  if and only if 
\begin{equation}\label{Rcompat}(\extd\tens\id-\id\wedge\nabla_E)\sigma_E(e\tens\extd a)=(\id\wedge\sigma)(\nabla_E e\tens \extd a),\quad\forall e\in E,\ a\in A.\end{equation}
In this case, applying this equation to $e.c\tens\extd a$ for $c\in A$ instead of
$e\tens\extd a$ gives
\begin{align*}
&\kern-20pt(\extd\tens\id-\id\wedge\nabla_E)\sigma_E(e.c \tens \extd a)   = 
(\id\wedge\sigma_E)(\nabla_E (e.c)\tens\extd a)\cr
   =&\ 
(\id\wedge\sigma_E)(\nabla_E e\tens c.\extd a)
+(\id\wedge\sigma_E)(\sigma_E(e\tens \extd c)\tens\extd a).
\end{align*}
Using extendability this is
\begin{eqnarray*}
(\extd\tens\id-\id\wedge\nabla_E)\sigma_E(e \tens c.\extd a)   \,=\,
(\id\wedge\sigma_E)(\nabla_E e\tens c.\extd a)
+\sigma_E(e\tens \extd c\wedge\extd a),
\end{eqnarray*}
or for $\xi=c.\extd a\in\Omega^1$,
\begin{eqnarray*}
(\extd\tens\id-\id\wedge\nabla_E)\sigma_E(e \tens \xi)   \,=\,
(\id\wedge\sigma_E)(\nabla_E e\tens \xi)
+\sigma_E(e\tens \extd \xi).
\end{eqnarray*}
From this $n=1$ case we prove the general equality by induction on $n$: Suppose that it is true for $n$, and
consider it for $n+1$. It is enough to check it on products $\xi\wedge\eta$, for
$\xi\in\Omega^1$ and $\eta\in\Omega^n$,
\begin{eqnarray} \label{yyiioo}
\nabla_E{}^{[n+1]}\sigma_E(e\tens\xi\wedge\eta) &=& \nabla_E{}^{[n+1]}(\wedge\tens\id)(\id\tens\sigma_E)
(\sigma_E(e\tens\xi)\tens \eta) \cr
&=& (\extd\wedge\sigma_E-\id\wedge\nabla_E{}^{[n]}\sigma_E)
(\sigma_E(e\tens\xi)\tens\eta)
\end{eqnarray}
by using extendability. Now by the $n=1$ case and more extendability,
\begin{align*}
&\kern-7pt(\extd\tens\id)\sigma_E(e\tens\xi) = (\id\wedge\nabla_E)\sigma_E(e\tens\xi)+
(\id\wedge\sigma_E)(\nabla_E e\tens\xi)+\sigma_E(e\tens\extd\xi),\cr
&\kern-7pt(\extd\wedge\sigma_E)(\sigma_E(e\tens\xi)\tens\eta) =
 (\id\wedge\sigma_E)((\id\wedge\nabla_E)\sigma_E(e\tens\xi)\tens\eta) \cr
&\quad +\, (\id\wedge\sigma_E)((\id\wedge\sigma_E)(\nabla_E e\tens\xi)\tens\eta)
+(\id\wedge\sigma_E)(\sigma_E(e\tens\extd\xi)\tens\eta) \cr
&= (\id\wedge\sigma_E)((\id\wedge\nabla_E)\sigma_E(e\tens\xi)\tens\eta) 
+ (\id\wedge\sigma_E)(\nabla_E e\tens\xi\wedge\eta)+\sigma_E(e\tens\extd\xi\wedge\eta) .
\end{align*}
Also the inductive hypothesis for $n$ and more extendability gives
\begin{align*}
&\kern-10pt (\id\wedge\nabla_E{}^{[n]}\sigma_E)
(\sigma_E(e\tens\xi)\tens\eta) \cr 
=&\ 
(\id\wedge((\id\wedge\sigma_E)(\nabla_E\tens\id)+\sigma_E(\id\tens\extd)))(\sigma_E(e\tens\xi)\tens\eta) \cr
=&\ (\id\wedge\id\wedge\sigma_E)
(\id\tens\nabla_E\tens\id)(\sigma_E(e\tens\xi)\tens\eta) 
 + (\id\wedge\sigma_E)(\sigma_E\tens\id)(e\tens\xi\tens\extd\eta) \cr
=&\ (\id\wedge\sigma_E)((\id\wedge\nabla_E)\sigma_E(e\tens\xi)\tens\eta)+\sigma_E(e\tens\xi\wedge
\extd\eta).
\end{align*}
Combining these in (\ref{yyiioo}) yields the required result,
\begin{eqnarray*}
\nabla_E{}^{[n+1]}\sigma_E(e\tens\xi\wedge\eta) = 
(\id\wedge\sigma_E)(\nabla_E e\tens\xi\wedge\eta)+\sigma_E(e\tens\extd\xi\wedge\eta)-
\sigma_E(e\tens\xi\wedge
\extd\eta).
\end{eqnarray*}
For the converse, the $n=1$ case of the equation here implies that $R_E$ is a right module map by condition (\ref{Rcompat}) since $\extd^2=0$.  Also in this case
\begin{align*}
&\kern-10pt(\id\wedge R_E)\sigma_E = (\id\wedge\nabla_E{}^{[1]})(\id\tens\nabla_E)\sigma_E\cr
=&\ (\id\wedge\nabla_E{}^{[1]})  (\extd\tens\id)\sigma_E- (\id\wedge\nabla^{[1]}) \nabla_E{}^{[1]}\sigma_E 
=  \nabla_E{}^{[2]} \nabla_E{}^{[1]}\sigma_E   \cr
=&\  \nabla_E{}^{[2]} (\id\wedge\sigma_E)(\nabla_E\tens\id)+ \nabla_E{}^{[2]} \sigma_E(\id\tens\extd)    \cr
=&\  \nabla_E{}^{[2]} (\id\wedge\sigma_E)(\nabla_E\tens\id)+ (\id\wedge\sigma_E)(\nabla_E\tens\extd)  
+ \sigma_E(\id\tens\extd)(\id\tens\extd)    \cr
=&\  \nabla_E{}^{[2]} (\id\wedge\sigma_E)(\nabla_E\tens\id)+ (\id\wedge\sigma_E)(\nabla_E\tens\extd)     \cr
=&\  (\extd\wedge\sigma_E)(\nabla_E\tens\id)- (\id\wedge\nabla_E{}^{[1]} \sigma_E)(\nabla_E\tens\id)
+ (\id\wedge\sigma_E)(\nabla_E\tens\extd)     \cr
=&\  (\extd\wedge\sigma_E)(\nabla_E\tens\id)- (\id\wedge(\id\wedge\sigma_E)(\nabla_E\tens\id))(\nabla_E\tens\id)\cr
& -\ (\id\wedge\sigma_E(\id\tens\extd))(\nabla_E\tens\id)
+ (\id\wedge\sigma_E)(\nabla_E\tens\extd)     \cr
=&\  (\extd\wedge\sigma_E)(\nabla_E\tens\id)- (\id\wedge(\id\wedge\sigma_E)(\nabla_E\tens\id))(\nabla_E\tens\id)\cr
=&\ (\id\wedge\sigma_E)((\extd\tens\id)\nabla_E\tens\id)-(\id\wedge\sigma_E)
((\id\wedge\nabla_E)\nabla_E\tens\id)\cr
=&\  (\id\wedge\sigma_E)(\nabla^{[1]} \nabla_E\tens\id)\,=\,  (\id\wedge\sigma_E)(R_E\tens\id).
\end{align*}
Expressing the stated condition in diagrammatic form

\medskip
\unitlength 0.56 mm
\begin{picture}(100,66)(-55,33)
\linethickness{0.3mm}
\put(10,85){\line(0,1){5}}
\linethickness{0.3mm}
\multiput(10,85)(0.12,-0.12){167}{\line(1,0){0.12}}
\linethickness{0.3mm}
\put(30,85){\line(0,1){5}}
\linethickness{0.3mm}
\multiput(21.25,77.5)(0.14,0.12){63}{\line(1,0){0.14}}
\linethickness{0.3mm}
\multiput(10,65)(0.12,0.14){63}{\line(0,1){0.14}}
\linethickness{0.3mm}
\multiput(39.99,55.5)(0.01,-0.5){1}{\line(0,-1){0.5}}
\multiput(39.95,56)(0.04,-0.5){1}{\line(0,-1){0.5}}
\multiput(39.89,56.49)(0.06,-0.49){1}{\line(0,-1){0.49}}
\multiput(39.8,56.98)(0.09,-0.49){1}{\line(0,-1){0.49}}
\multiput(39.69,57.47)(0.11,-0.49){1}{\line(0,-1){0.49}}
\multiput(39.56,57.95)(0.14,-0.48){1}{\line(0,-1){0.48}}
\multiput(39.4,58.42)(0.16,-0.47){1}{\line(0,-1){0.47}}
\multiput(39.21,58.88)(0.09,-0.23){2}{\line(0,-1){0.23}}
\multiput(39.01,59.34)(0.1,-0.23){2}{\line(0,-1){0.23}}
\multiput(38.78,59.78)(0.11,-0.22){2}{\line(0,-1){0.22}}
\multiput(38.53,60.21)(0.12,-0.22){2}{\line(0,-1){0.22}}
\multiput(38.26,60.63)(0.14,-0.21){2}{\line(0,-1){0.21}}
\multiput(37.97,61.04)(0.15,-0.2){2}{\line(0,-1){0.2}}
\multiput(37.66,61.43)(0.1,-0.13){3}{\line(0,-1){0.13}}
\multiput(37.33,61.8)(0.11,-0.12){3}{\line(0,-1){0.12}}
\multiput(36.98,62.16)(0.12,-0.12){3}{\line(0,-1){0.12}}
\multiput(36.62,62.5)(0.12,-0.11){3}{\line(1,0){0.12}}
\multiput(36.23,62.82)(0.13,-0.11){3}{\line(1,0){0.13}}
\multiput(35.84,63.12)(0.13,-0.1){3}{\line(1,0){0.13}}
\multiput(35.43,63.4)(0.21,-0.14){2}{\line(1,0){0.21}}
\multiput(35,63.66)(0.21,-0.13){2}{\line(1,0){0.21}}
\multiput(34.56,63.9)(0.22,-0.12){2}{\line(1,0){0.22}}
\multiput(34.11,64.12)(0.22,-0.11){2}{\line(1,0){0.22}}
\multiput(33.65,64.31)(0.23,-0.1){2}{\line(1,0){0.23}}
\multiput(33.18,64.48)(0.47,-0.17){1}{\line(1,0){0.47}}
\multiput(32.71,64.63)(0.48,-0.15){1}{\line(1,0){0.48}}
\multiput(32.23,64.75)(0.48,-0.12){1}{\line(1,0){0.48}}
\multiput(31.74,64.85)(0.49,-0.1){1}{\line(1,0){0.49}}
\multiput(31.24,64.92)(0.49,-0.07){1}{\line(1,0){0.49}}
\multiput(30.75,64.97)(0.5,-0.05){1}{\line(1,0){0.5}}
\multiput(30.25,65)(0.5,-0.02){1}{\line(1,0){0.5}}
\put(29.75,65){\line(1,0){0.5}}
\multiput(29.25,64.97)(0.5,0.02){1}{\line(1,0){0.5}}
\multiput(28.76,64.92)(0.5,0.05){1}{\line(1,0){0.5}}
\multiput(28.26,64.85)(0.49,0.07){1}{\line(1,0){0.49}}
\multiput(27.77,64.75)(0.49,0.1){1}{\line(1,0){0.49}}
\multiput(27.29,64.63)(0.48,0.12){1}{\line(1,0){0.48}}
\multiput(26.82,64.48)(0.48,0.15){1}{\line(1,0){0.48}}
\multiput(26.35,64.31)(0.47,0.17){1}{\line(1,0){0.47}}
\multiput(25.89,64.12)(0.23,0.1){2}{\line(1,0){0.23}}
\multiput(25.44,63.9)(0.22,0.11){2}{\line(1,0){0.22}}
\multiput(25,63.66)(0.22,0.12){2}{\line(1,0){0.22}}
\multiput(24.57,63.4)(0.21,0.13){2}{\line(1,0){0.21}}
\multiput(24.16,63.12)(0.21,0.14){2}{\line(1,0){0.21}}
\multiput(23.77,62.82)(0.13,0.1){3}{\line(1,0){0.13}}
\multiput(23.38,62.5)(0.13,0.11){3}{\line(1,0){0.13}}
\multiput(23.02,62.16)(0.12,0.11){3}{\line(1,0){0.12}}
\multiput(22.67,61.8)(0.12,0.12){3}{\line(0,1){0.12}}
\multiput(22.34,61.43)(0.11,0.12){3}{\line(0,1){0.12}}
\multiput(22.03,61.04)(0.1,0.13){3}{\line(0,1){0.13}}
\multiput(21.74,60.63)(0.15,0.2){2}{\line(0,1){0.2}}
\multiput(21.47,60.21)(0.14,0.21){2}{\line(0,1){0.21}}
\multiput(21.22,59.78)(0.12,0.22){2}{\line(0,1){0.22}}
\multiput(20.99,59.34)(0.11,0.22){2}{\line(0,1){0.22}}
\multiput(20.79,58.88)(0.1,0.23){2}{\line(0,1){0.23}}
\multiput(20.6,58.42)(0.09,0.23){2}{\line(0,1){0.23}}
\multiput(20.44,57.95)(0.16,0.47){1}{\line(0,1){0.47}}
\multiput(20.31,57.47)(0.14,0.48){1}{\line(0,1){0.48}}
\multiput(20.2,56.98)(0.11,0.49){1}{\line(0,1){0.49}}
\multiput(20.11,56.49)(0.09,0.49){1}{\line(0,1){0.49}}
\multiput(20.05,56)(0.06,0.49){1}{\line(0,1){0.49}}
\multiput(20.01,55.5)(0.04,0.5){1}{\line(0,1){0.5}}
\multiput(20,55)(0.01,0.5){1}{\line(0,1){0.5}}

\linethickness{0.3mm}
\multiput(10,55)(0.03,-0.51){1}{\line(0,-1){0.51}}
\multiput(10.03,54.49)(0.08,-0.5){1}{\line(0,-1){0.5}}
\multiput(10.1,53.99)(0.13,-0.49){1}{\line(0,-1){0.49}}
\multiput(10.23,53.5)(0.18,-0.47){1}{\line(0,-1){0.47}}
\multiput(10.41,53.03)(0.11,-0.23){2}{\line(0,-1){0.23}}
\multiput(10.63,52.57)(0.13,-0.21){2}{\line(0,-1){0.21}}
\multiput(10.9,52.14)(0.1,-0.13){3}{\line(0,-1){0.13}}
\multiput(11.21,51.74)(0.12,-0.12){3}{\line(0,-1){0.12}}
\multiput(11.56,51.38)(0.13,-0.11){3}{\line(1,0){0.13}}
\multiput(11.94,51.05)(0.21,-0.14){2}{\line(1,0){0.21}}
\multiput(12.36,50.76)(0.22,-0.12){2}{\line(1,0){0.22}}
\multiput(12.8,50.51)(0.23,-0.1){2}{\line(1,0){0.23}}
\multiput(13.26,50.31)(0.48,-0.15){1}{\line(1,0){0.48}}
\multiput(13.75,50.16)(0.5,-0.1){1}{\line(1,0){0.5}}
\multiput(14.24,50.06)(0.5,-0.05){1}{\line(1,0){0.5}}
\put(14.75,50.01){\line(1,0){0.51}}
\multiput(15.25,50.01)(0.5,0.05){1}{\line(1,0){0.5}}
\multiput(15.76,50.06)(0.5,0.1){1}{\line(1,0){0.5}}
\multiput(16.25,50.16)(0.48,0.15){1}{\line(1,0){0.48}}
\multiput(16.74,50.31)(0.23,0.1){2}{\line(1,0){0.23}}
\multiput(17.2,50.51)(0.22,0.12){2}{\line(1,0){0.22}}
\multiput(17.64,50.76)(0.21,0.14){2}{\line(1,0){0.21}}
\multiput(18.06,51.05)(0.13,0.11){3}{\line(1,0){0.13}}
\multiput(18.44,51.38)(0.12,0.12){3}{\line(0,1){0.12}}
\multiput(18.79,51.74)(0.1,0.13){3}{\line(0,1){0.13}}
\multiput(19.1,52.14)(0.13,0.21){2}{\line(0,1){0.21}}
\multiput(19.37,52.57)(0.11,0.23){2}{\line(0,1){0.23}}
\multiput(19.59,53.03)(0.18,0.47){1}{\line(0,1){0.47}}
\multiput(19.77,53.5)(0.13,0.49){1}{\line(0,1){0.49}}
\multiput(19.9,53.99)(0.08,0.5){1}{\line(0,1){0.5}}
\multiput(19.97,54.49)(0.03,0.51){1}{\line(0,1){0.51}}

\linethickness{0.3mm}
\put(15,40){\line(0,1){10}}
\linethickness{0.3mm}
\put(10,55){\line(0,1){10}}
\linethickness{0.3mm}
\put(40,40){\line(0,1){15}}
\linethickness{0.3mm}
\put(70,85){\line(0,1){5}}
\linethickness{0.3mm}
\multiput(79.99,75.5)(0.01,-0.5){1}{\line(0,-1){0.5}}
\multiput(79.95,76)(0.04,-0.5){1}{\line(0,-1){0.5}}
\multiput(79.89,76.49)(0.06,-0.49){1}{\line(0,-1){0.49}}
\multiput(79.8,76.98)(0.09,-0.49){1}{\line(0,-1){0.49}}
\multiput(79.69,77.47)(0.11,-0.49){1}{\line(0,-1){0.49}}
\multiput(79.56,77.95)(0.14,-0.48){1}{\line(0,-1){0.48}}
\multiput(79.4,78.42)(0.16,-0.47){1}{\line(0,-1){0.47}}
\multiput(79.21,78.88)(0.09,-0.23){2}{\line(0,-1){0.23}}
\multiput(79.01,79.34)(0.1,-0.23){2}{\line(0,-1){0.23}}
\multiput(78.78,79.78)(0.11,-0.22){2}{\line(0,-1){0.22}}
\multiput(78.53,80.21)(0.12,-0.22){2}{\line(0,-1){0.22}}
\multiput(78.26,80.63)(0.14,-0.21){2}{\line(0,-1){0.21}}
\multiput(77.97,81.04)(0.15,-0.2){2}{\line(0,-1){0.2}}
\multiput(77.66,81.43)(0.1,-0.13){3}{\line(0,-1){0.13}}
\multiput(77.33,81.8)(0.11,-0.12){3}{\line(0,-1){0.12}}
\multiput(76.98,82.16)(0.12,-0.12){3}{\line(0,-1){0.12}}
\multiput(76.62,82.5)(0.12,-0.11){3}{\line(1,0){0.12}}
\multiput(76.23,82.82)(0.13,-0.11){3}{\line(1,0){0.13}}
\multiput(75.84,83.12)(0.13,-0.1){3}{\line(1,0){0.13}}
\multiput(75.43,83.4)(0.21,-0.14){2}{\line(1,0){0.21}}
\multiput(75,83.66)(0.21,-0.13){2}{\line(1,0){0.21}}
\multiput(74.56,83.9)(0.22,-0.12){2}{\line(1,0){0.22}}
\multiput(74.11,84.12)(0.22,-0.11){2}{\line(1,0){0.22}}
\multiput(73.65,84.31)(0.23,-0.1){2}{\line(1,0){0.23}}
\multiput(73.18,84.48)(0.47,-0.17){1}{\line(1,0){0.47}}
\multiput(72.71,84.63)(0.48,-0.15){1}{\line(1,0){0.48}}
\multiput(72.23,84.75)(0.48,-0.12){1}{\line(1,0){0.48}}
\multiput(71.74,84.85)(0.49,-0.1){1}{\line(1,0){0.49}}
\multiput(71.24,84.92)(0.49,-0.07){1}{\line(1,0){0.49}}
\multiput(70.75,84.97)(0.5,-0.05){1}{\line(1,0){0.5}}
\multiput(70.25,85)(0.5,-0.02){1}{\line(1,0){0.5}}
\put(69.75,85){\line(1,0){0.5}}
\multiput(69.25,84.97)(0.5,0.02){1}{\line(1,0){0.5}}
\multiput(68.76,84.92)(0.5,0.05){1}{\line(1,0){0.5}}
\multiput(68.26,84.85)(0.49,0.07){1}{\line(1,0){0.49}}
\multiput(67.77,84.75)(0.49,0.1){1}{\line(1,0){0.49}}
\multiput(67.29,84.63)(0.48,0.12){1}{\line(1,0){0.48}}
\multiput(66.82,84.48)(0.48,0.15){1}{\line(1,0){0.48}}
\multiput(66.35,84.31)(0.47,0.17){1}{\line(1,0){0.47}}
\multiput(65.89,84.12)(0.23,0.1){2}{\line(1,0){0.23}}
\multiput(65.44,83.9)(0.22,0.11){2}{\line(1,0){0.22}}
\multiput(65,83.66)(0.22,0.12){2}{\line(1,0){0.22}}
\multiput(64.57,83.4)(0.21,0.13){2}{\line(1,0){0.21}}
\multiput(64.16,83.12)(0.21,0.14){2}{\line(1,0){0.21}}
\multiput(63.77,82.82)(0.13,0.1){3}{\line(1,0){0.13}}
\multiput(63.38,82.5)(0.13,0.11){3}{\line(1,0){0.13}}
\multiput(63.02,82.16)(0.12,0.11){3}{\line(1,0){0.12}}
\multiput(62.67,81.8)(0.12,0.12){3}{\line(0,1){0.12}}
\multiput(62.34,81.43)(0.11,0.12){3}{\line(0,1){0.12}}
\multiput(62.03,81.04)(0.1,0.13){3}{\line(0,1){0.13}}
\multiput(61.74,80.63)(0.15,0.2){2}{\line(0,1){0.2}}
\multiput(61.47,80.21)(0.14,0.21){2}{\line(0,1){0.21}}
\multiput(61.22,79.78)(0.12,0.22){2}{\line(0,1){0.22}}
\multiput(60.99,79.34)(0.11,0.22){2}{\line(0,1){0.22}}
\multiput(60.79,78.88)(0.1,0.23){2}{\line(0,1){0.23}}
\multiput(60.6,78.42)(0.09,0.23){2}{\line(0,1){0.23}}
\multiput(60.44,77.95)(0.16,0.47){1}{\line(0,1){0.47}}
\multiput(60.31,77.47)(0.14,0.48){1}{\line(0,1){0.48}}
\multiput(60.2,76.98)(0.11,0.49){1}{\line(0,1){0.49}}
\multiput(60.11,76.49)(0.09,0.49){1}{\line(0,1){0.49}}
\multiput(60.05,76)(0.06,0.49){1}{\line(0,1){0.49}}
\multiput(60.01,75.5)(0.04,0.5){1}{\line(0,1){0.5}}
\multiput(60,75)(0.01,0.5){1}{\line(0,1){0.5}}

\linethickness{0.3mm}
\multiput(80,75)(0.12,-0.12){167}{\line(1,0){0.12}}
\linethickness{0.3mm}
\put(100,80){\line(0,1){10}}
\linethickness{0.3mm}
\multiput(90,67.5)(0.12,0.15){83}{\line(0,1){0.15}}
\linethickness{0.3mm}
\multiput(60,55)(0.01,-0.5){1}{\line(0,-1){0.5}}
\multiput(60.01,54.5)(0.04,-0.5){1}{\line(0,-1){0.5}}
\multiput(60.05,54)(0.06,-0.49){1}{\line(0,-1){0.49}}
\multiput(60.11,53.51)(0.09,-0.49){1}{\line(0,-1){0.49}}
\multiput(60.2,53.02)(0.11,-0.49){1}{\line(0,-1){0.49}}
\multiput(60.31,52.53)(0.14,-0.48){1}{\line(0,-1){0.48}}
\multiput(60.44,52.05)(0.16,-0.47){1}{\line(0,-1){0.47}}
\multiput(60.6,51.58)(0.09,-0.23){2}{\line(0,-1){0.23}}
\multiput(60.79,51.12)(0.1,-0.23){2}{\line(0,-1){0.23}}
\multiput(60.99,50.66)(0.11,-0.22){2}{\line(0,-1){0.22}}
\multiput(61.22,50.22)(0.12,-0.22){2}{\line(0,-1){0.22}}
\multiput(61.47,49.79)(0.14,-0.21){2}{\line(0,-1){0.21}}
\multiput(61.74,49.37)(0.15,-0.2){2}{\line(0,-1){0.2}}
\multiput(62.03,48.96)(0.1,-0.13){3}{\line(0,-1){0.13}}
\multiput(62.34,48.57)(0.11,-0.12){3}{\line(0,-1){0.12}}
\multiput(62.67,48.2)(0.12,-0.12){3}{\line(0,-1){0.12}}
\multiput(63.02,47.84)(0.12,-0.11){3}{\line(1,0){0.12}}
\multiput(63.38,47.5)(0.13,-0.11){3}{\line(1,0){0.13}}
\multiput(63.77,47.18)(0.13,-0.1){3}{\line(1,0){0.13}}
\multiput(64.16,46.88)(0.21,-0.14){2}{\line(1,0){0.21}}
\multiput(64.57,46.6)(0.21,-0.13){2}{\line(1,0){0.21}}
\multiput(65,46.34)(0.22,-0.12){2}{\line(1,0){0.22}}
\multiput(65.44,46.1)(0.22,-0.11){2}{\line(1,0){0.22}}
\multiput(65.89,45.88)(0.23,-0.1){2}{\line(1,0){0.23}}
\multiput(66.35,45.69)(0.47,-0.17){1}{\line(1,0){0.47}}
\multiput(66.82,45.52)(0.48,-0.15){1}{\line(1,0){0.48}}
\multiput(67.29,45.37)(0.48,-0.12){1}{\line(1,0){0.48}}
\multiput(67.77,45.25)(0.49,-0.1){1}{\line(1,0){0.49}}
\multiput(68.26,45.15)(0.49,-0.07){1}{\line(1,0){0.49}}
\multiput(68.76,45.08)(0.5,-0.05){1}{\line(1,0){0.5}}
\multiput(69.25,45.03)(0.5,-0.02){1}{\line(1,0){0.5}}
\put(69.75,45){\line(1,0){0.5}}
\multiput(70.25,45)(0.5,0.02){1}{\line(1,0){0.5}}
\multiput(70.75,45.03)(0.5,0.05){1}{\line(1,0){0.5}}
\multiput(71.24,45.08)(0.49,0.07){1}{\line(1,0){0.49}}
\multiput(71.74,45.15)(0.49,0.1){1}{\line(1,0){0.49}}
\multiput(72.23,45.25)(0.48,0.12){1}{\line(1,0){0.48}}
\multiput(72.71,45.37)(0.48,0.15){1}{\line(1,0){0.48}}
\multiput(73.18,45.52)(0.47,0.17){1}{\line(1,0){0.47}}
\multiput(73.65,45.69)(0.23,0.1){2}{\line(1,0){0.23}}
\multiput(74.11,45.88)(0.22,0.11){2}{\line(1,0){0.22}}
\multiput(74.56,46.1)(0.22,0.12){2}{\line(1,0){0.22}}
\multiput(75,46.34)(0.21,0.13){2}{\line(1,0){0.21}}
\multiput(75.43,46.6)(0.21,0.14){2}{\line(1,0){0.21}}
\multiput(75.84,46.88)(0.13,0.1){3}{\line(1,0){0.13}}
\multiput(76.23,47.18)(0.13,0.11){3}{\line(1,0){0.13}}
\multiput(76.62,47.5)(0.12,0.11){3}{\line(1,0){0.12}}
\multiput(76.98,47.84)(0.12,0.12){3}{\line(0,1){0.12}}
\multiput(77.33,48.2)(0.11,0.12){3}{\line(0,1){0.12}}
\multiput(77.66,48.57)(0.1,0.13){3}{\line(0,1){0.13}}
\multiput(77.97,48.96)(0.15,0.2){2}{\line(0,1){0.2}}
\multiput(78.26,49.37)(0.14,0.21){2}{\line(0,1){0.21}}
\multiput(78.53,49.79)(0.12,0.22){2}{\line(0,1){0.22}}
\multiput(78.78,50.22)(0.11,0.22){2}{\line(0,1){0.22}}
\multiput(79.01,50.66)(0.1,0.23){2}{\line(0,1){0.23}}
\multiput(79.21,51.12)(0.09,0.23){2}{\line(0,1){0.23}}
\multiput(79.4,51.58)(0.16,0.47){1}{\line(0,1){0.47}}
\multiput(79.56,52.05)(0.14,0.48){1}{\line(0,1){0.48}}
\multiput(79.69,52.53)(0.11,0.49){1}{\line(0,1){0.49}}
\multiput(79.8,53.02)(0.09,0.49){1}{\line(0,1){0.49}}
\multiput(79.89,53.51)(0.06,0.49){1}{\line(0,1){0.49}}
\multiput(79.95,54)(0.04,0.5){1}{\line(0,1){0.5}}
\multiput(79.99,54.5)(0.01,0.5){1}{\line(0,1){0.5}}

\linethickness{0.3mm}
\put(70,40){\line(0,1){5}}
\linethickness{0.3mm}
\multiput(80,55)(0.12,0.14){63}{\line(0,1){0.14}}
\linethickness{0.3mm}
\put(60,55){\line(0,1){20}}
\linethickness{0.3mm}
\put(100,40){\line(0,1){15}}
\put(10,93){\makebox(0,0)[cc]{$E$}}

\put(30,93){\makebox(0,0)[cc]{$\Omega^n$}}

\put(100,93){\makebox(0,0)[cc]{$\Omega^n$}}

\put(70,93){\makebox(0,0)[cc]{$E$}}

\put(70,79){\makebox(0,0)[cc]{$R_E$}}

\put(31,58){\makebox(0,0)[cc]{$R_E$}}

\put(96,66.25){\makebox(0,0)[cc]{$\sigma_E$}}

\put(27,75){\makebox(0,0)[cc]{$\sigma_E$}}

\put(18.75,47){\makebox(0,0)[cc]{$\wedge$}}

\put(70,48.75){\makebox(0,0)[cc]{$\wedge$}}

\put(15,36){\makebox(0,0)[cc]{$\Omega^{n+2}$}}

\put(70,36){\makebox(0,0)[cc]{$\Omega^{n+2}$}}

\put(40,36){\makebox(0,0)[cc]{$E$}}

\put(100,36){\makebox(0,0)[cc]{$E$}}

\put(51,65){\makebox(0,0)[cc]{$=$}}

\end{picture}

\noindent
is the easiest way to extended to other $n$ by induction, a  tedious but straightforward exercise in associativity of wedge products. 
\eproof

\medskip
 We give a class of examples where extendability of the bimodule connections is automatic once the curvature is assumed to be a bimodule map. Recall that the maximal prolongation differential calculus extends a given first order differential calculus to all orders by imposing quadratic relations which are precisely $\extd a\wedge\extd b+\extd r\wedge\extd s=0\in\Omega^2$ for every relation
$a.\extd b=\extd r.s\in\Omega^1$ (summation on both sides implicit, for some $a,b,r,s\in A$).

\begin{lemma} \label{dinnu}
Suppose that $A$ is given the maximal prolongation differential calculus for some first order differential calculus and that $(E,\nabla_E,\sigma_E)$ is a left bimodule connection whose curvature $R_E$ is a right module map. Then $(E,\nabla_E,\sigma_E)$ is an extendable bimodule connection. 
\end{lemma}\proof Suppose that $a.\extd b=\extd r.s\in\Omega^1$ (summation implicit, for some $a,b,r,s\in A$). We also set $\sigma_E(e\tens\extd r)=\xi\tens f$ (summation understood) as a shorthand. Then
\begin{align*}
&\kern-20pt(\id\wedge\sigma_E)(\sigma_E(e\tens\extd r)\tens\extd s) 
 = \xi\wedge (\nabla_E(f.s)-(\nabla_Ef).s) \cr
=&\ \extd\xi\tens f.s -  (\extd\tens\id-\id\wedge\nabla_E)(\xi\tens f.s) - \xi\wedge (\nabla_E f).s \cr
=&\ (\extd\tens\id-\id\wedge\nabla_E)(\xi\tens f).s -  (\extd\tens\id-\id\wedge\nabla_E)(\xi\tens f.s) \cr
=&\ (\extd\tens\id-\id\wedge\nabla_E)\sigma_E(e\tens\extd r).s -  (\extd\tens\id-\id\wedge\nabla_E)\sigma_E(e\tens\extd r.s)\cr
=&\ (\extd\tens\id-\id\wedge\nabla_E)\sigma_E(e\tens\extd r).s -  (\extd\tens\id-\id\wedge\nabla_E)\sigma_E(e\tens a.\extd b)\cr
=&\ (\extd\tens\id-\id\wedge\nabla_E)\sigma_E(e\tens\extd r).s -  (\extd\tens\id-\id\wedge\nabla_E)\sigma_E(e.a\tens\extd b).
\end{align*}
If we use the expression for $R_E(e.a)-R_E(e).a$ from the proof of Lemma~\ref{mor2help} then
\begin{align*}
&\kern-20pt(\id\wedge\sigma_E)(\sigma_E(e\tens\extd r)\tens\extd s) \cr
=& (\id\wedge\sigma_E)  (\nabla_E e\tens \extd r.s  -\nabla_E(e.a) \tens\extd b) + R_E(er)s-R_E(e)rs-R_E(eab)+R_E(ea)b  \cr
=&(\id\wedge\sigma_E)  (\nabla_E e\tens a.\extd b  -\nabla_E(e.a) \tens\extd b) + R_E(er)s-R_E(e)rs-R_E(eab)+R_E(ea)b  \cr
=&(\id\wedge\sigma_E) ( (\nabla_E( e).a  -\nabla_E(e.a)) \tens\extd b)  + R_E(er)s-R_E(e)rs-R_E(eab)+R_E(ea)b \cr
=&-\,(\id\wedge\sigma_E) (\sigma_E( e\tens\extd a ) \tens\extd b) + R_E(er)s-R_E(e)rs-R_E(eab)+R_E(ea)b  .
\end{align*}
Hence
\begin{eqnarray*}
(\id\wedge\sigma_E)(\sigma_E\tens\id)(e\tens(\extd a\tens\extd b+\extd r\tens\extd s))= R_E(er)s-R_E(e)rs-R_E(eab)+R_E(ea)b.
\end{eqnarray*}
In particular, if $R_E$ is a right module map then the LHS vanishes for all $e$ and all $a,b,r,s$ (sum of such) obeying the condition stated. The relations for the wedge product are given by 
$\extd a\tens\extd b+\extd r\tens\extd s$  whenever
$a.\extd b=\extd r.s\in\Omega^1$ (summation implicit). Hence if $R_E$ is a right module map, for 
$(\id\tens\sigma_E)(\sigma_E\tens\id):E\tens_A\Omega^1
\tens_A\Omega^1\to \Omega^1\tens_A\Omega^1\tens_AE$
then we find 
\begin{equation}\label{maxprocalc}
(\id\tens\sigma_E)(\sigma_E\tens\id)(E\otimes_A \ker\wedge)
\subset \ker\wedge \otimes_A E.
\end{equation}
As $\Omega^2$ is $\Omega^1 \otimes_A   \Omega^1$ quotiented by $\ker\wedge$, we find that $(\id\tens\sigma_E)(\sigma_E\tens\id)$ induces a map 
$\sigma_E:E\tens_A\Omega^2\to \Omega^2 \tens_A E$ so we have proven extendability in degree 2. The higher extendability in the case of the maximal prolongation is then automatic. Thus for  $\Omega^3$ consider the map
\[
(\id^{\tens 2}\tens \sigma_E)(\id\tens\sigma_E\tens\id)(\sigma_E\tens\id^{\tens 2}):E\otimes_A\Omega^1\otimes_A\Omega^1\otimes_A\Omega^1 \to 
\Omega^1\otimes_A\Omega^1\otimes_A\Omega^1\otimes_A E
\]
By (\ref{maxprocalc}), this sends $E\tens_A \ker\wedge\tens_A \Omega^1$ to $\ker\wedge\tens_A E\tens_A \Omega^1$ and sends
$E\tens_A \Omega^1\tens_A\ker\wedge$ to $\Omega^1\tens_A E\tens_A\ker\wedge$. Then by definition of maximal prolongation there is a well defined map from $E\tens_A\Omega^3$ to $\Omega^3\tens_A E$. Similary for higher degree. \eproof

We can use extendability to give a simple statement about the curvature of a tensor product:

\begin{lemma} \label{fourteen} 
For a left $A$-module $F$ and an $A$-bimodule $E$, suppose that $(E,\nabla_E,\sigma_E)$ is an extendable bimodule connection whose curvature is a right module map and that
 $(F,\nabla_F)$ is a left connection. Then
\[
R_{E\tens F}=R_E\tens\id+(\sigma_E\tens\id)(\id\tens R_F): E\otimes_A F\to \Omega^2\otimes_A E\otimes_A F
\]
\end{lemma}
\proof Using Lemma~\ref{mor2help}:
\begin{align*}
&\kern-20pt \nabla_{E\tens F}{}^{[1]}\,\nabla^{\phantom{[1]}}_{E\tens F} = \nabla_{E\tens F}{}^{[1]}\,
(\nabla_E\tens\id+(\sigma_E\tens\id)(\id\tens\nabla_F)) \cr
=&\ \nabla_{E}{}^{[1]}\nabla_E\tens\id-(\id\wedge\sigma_E\tens\id)(\nabla_E\tens\nabla_F) \cr
& +\, (\nabla_E{}^{[1]}\,\sigma_E\tens\id)(\id\tens\nabla_F)
-(\id\wedge\sigma_E\tens\id)(\sigma_E\tens\nabla_F)(\id\tens\nabla_F) \cr
=&\ \nabla_{E}{}^{[1]}\nabla_E\tens\id+(\sigma_E\tens\id)(\id\tens\extd\tens\id)
(\id\tens\nabla_F) \cr
& -\, (\id\wedge\sigma_E\tens\id)(\sigma_E\tens\nabla_F)(\id\tens\nabla_F) \cr
=&\ \nabla_{E}{}^{[1]}\nabla_E\tens\id+(\sigma_E\tens\id)(\id\tens\extd\tens\id)
(\id\tens\nabla_F) \cr
& -\, (\sigma_E\tens\id)(\id\tens(\id\wedge\nabla_F)\nabla_F)
\end{align*}
which gives the required formula
\eproof

\subsection{Construction of  ${}_A\mathcal{G}_A$}\label{secAGA2}
Consider two $A$-bimodules with left bimodule connections $(E,\nabla_E,\sigma_E)$ and $(F,\nabla_F,\sigma_F)$, and suppose that $\phi:E\to F$ is a bimodule map. We know that  $\doublenabla(\phi):E\to \Omega^1\tens_A F$ is a left module map, and we now check the right action:
\begin{align} \label{uytr}
\doublenabla(\phi )(e.a) =&\ \nabla_F\,\phi (e.a) - (\id\tens \phi )\,\nabla_E(e.a) 
= \nabla_F( (\phi e).a) - (\id\tens \phi )\,\nabla_E(e.a) \cr
=&(\nabla_F\phi e).a- ((\id\tens \phi )\,\nabla_Ee).a+\sigma_F(\phi (e)\tens \extd a) - \  (\id\tens \phi )\,\sigma_E(e\tens \extd a)\cr
=& (\doublenabla(\phi)(e)).a +\sigma_F(\phi e\tens \extd a)  - \  (\id\tens \phi )\,\sigma_E(e\tens \extd a)
\end{align} 
for $a\in A$ and $e\in E$. So $\doublenabla(\phi)$ is a bimodule map if and only if  $(\id\tens \phi )\,\sigma_E=\sigma_F\,(\phi \tens\id)$.
If we have $0$-morphisms being bimodule maps it is reasonable to expect that their derivatives, the $1$-morphisms, should also be bimodule maps. Thus we assume the additional condition $(\id\tens \phi )\,\sigma_E=\sigma_F\,(\phi \tens\id)$ for $0$-morphisms in our construction of a bimodule DG category. 
 The required generalisation for morphisms $\psi:E\to \Omega^n\tens_A F$ of grade $n$ is bimodule maps obeying the condition
\begin{equation}\label{AGAmorn} 
(\id\wedge\psi)\sigma_E=(-1)^n(\id\wedge\sigma_F)(\phi\tens\id):E\otimes_A\Omega^1\to \Omega^{n+1}\otimes_A F
\end{equation}
and Lemma~\ref{agadiff} will say that  $\doublenabla(\psi)$  is then an $n+1$-morphism. This gives the morphisms in our bimodule category 
${}_A\CG_A$. Remembering how important the curvature in the theory is, it is reasonable to restrict to objects whose curvature is a bimodule map. Also as we wish to take tensor products, it is reasonable to restrict to extendable bimodule connections as then we have the result of Lemma~\ref{fourteen} on the curvature of tensor products. (In fact the technical result Lemma~\ref{mor2help} will prove useful in numerous places.)
 We summarise our proposed DG category by the following table:
 
\medskip
\noindent
{\renewcommand{\arraystretch}{1}
\begin{tabular}{c|c|c}Name & Objects &   $\Mor_n((E,\nabla_E,\sigma_E),(F,\nabla_F,\sigma_F))$ \\
\hline
${}_A\mathcal{G}_{A}$\index{category!${}_A\mathcal{G}_{A}$}  & $(E,\nabla_E,\sigma_E)$ &   $\phi:E\to\Omega^n \otimes_A F$    \\
 & extendable left $A-A$-bimodule connections   & a bimodule map obeying (\ref{AGAmorn})  \\
 & with $R_E$ a bimodule map  &  
\end{tabular}}

\medskip
It is useful to note that by extendability and a use of induction a bimodule map $\phi$ obeying (\ref{AGAmorn}) also obeys
the following for all $n,m\ge 0$:
 \begin{align} \label{severalsign}
(\id\wedge\phi)\sigma_E=(-1)^{nm}  (\id\wedge\sigma_F)(\phi\tens\id) : E\otimes_A\Omega^m \to \Omega^{n+m}  \otimes_A F
\end{align}

\begin{theorem} \label{baacatdef} 
If we use the composition and derivative of morphisms as for ${}_A\mathcal{G}$, then
${}_A\mathcal{G}_{A}$ as summarised above is a DG category and the conclusions of Theorem~\ref{the-bcat} hold but now for ${}_A\CG_A$. 
\end{theorem}
\proof We begin by checking that ${}_A\mathcal{G}_{A}$ is a category, i.e. that we can compose morphisms. This is shown using diagramatic notation  to check that (\ref{AGAmorn}) holds for the composition of $\phi\in\Mor_n(E,F)$ and $\psi\in\Mor_m(F,G)$:

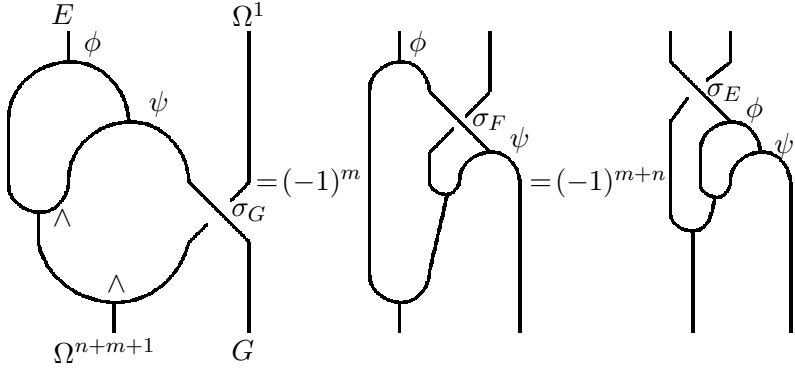
\unitlength 0.8 mm


\noindent
It is clear that  $\phi\circ\psi$ is a bimodule map if $\phi,\psi$ are. Lemma~\ref{mor2help} showed that the curvature is a morphism and Lemma~\ref{agadiff} will show that the differential of a morphism is a morphism. The rest follows just as for ${}_A\mathcal{G}$. \eproof

 To conclude the proof of Theorem~\ref{baacatdef} we need the following:

\begin{lemma} \label{agadiff}
If $\phi\in \mathrm{Mor}_n(E,F)$ in ${}_A\mathcal{G}_A$,
 then $\doublenabla(\phi)\in \mathrm{Mor}_{n+1}(E,F)$. 
\end{lemma}
\proof First we check that $\doublenabla(\phi)$ is a right module map. For all $e\in E$ and $a\in A$ and setting $\phi e=\xi\tens f$ (summation understood), 
\begin{align*}
\doublenabla(\phi)(e.a) =&\nabla^{[n]}_F(\phi(e.a))-(\id\wedge\phi)\nabla_E(e.a) \cr
=&
 \nabla^{[n]}_F(\xi\tens f.a)-(\id\wedge\phi)(\nabla_E e).a
-(\id\wedge\phi)\sigma_E(e\tens\extd a)\cr
=&\
 \extd\xi\tens f.a+ (-1)^n\ \xi\wedge \nabla_F(f.a)-(\id\wedge\phi)(\nabla_Ee).a
-(\id\wedge\phi)\sigma_E(e\tens\extd a)\cr
=&\
\doublenabla(\phi)(e).a+ (-1)^n\ \xi\wedge \sigma_F(f\tens\extd a)
-(\id\wedge\phi)\sigma_E(e\tens\extd a)  =
\doublenabla(\phi)(e).a
\end{align*}
where the last equality is the assumption (\ref{AGAmorn}). Now check that $\doublenabla(\phi)$ itself obeys the condition (\ref{AGAmorn}) in degree $n+1$. Here
\begin{eqnarray*}
(\id\wedge\doublenabla(\phi))\sigma_E = (\id\wedge\nabla^{[n]}_F\phi)\sigma_E
-(\id\wedge\phi)(\id \wedge\nabla_E)\sigma_E: E\otimes_A\Omega^1 \to \Omega^{n+2}  \otimes_A F.
\end{eqnarray*}
From the definition of $\nabla_F{}^{[n+1]}$,
\begin{eqnarray*}
\nabla_F^{[n+1]}(\id\wedge\phi)\sigma_E = (\extd \wedge\phi)\sigma_E-
 (\id\wedge\nabla^{[n]}_F\phi)\sigma_E,
\end{eqnarray*}
so using Lemma~\ref{mor2help},
\begin{align*}
&\kern-20pt(\id\wedge\doublenabla(\phi))\sigma_E = (\extd \wedge\phi)\sigma_E-\nabla_F^{[n+1]}(\id\wedge\phi)\sigma_E
-(\id\wedge\phi)(\id \wedge\nabla_E)\sigma_E   \cr
=&\ (\id\wedge\phi)\nabla_E{}^{[1]}\sigma_E
-(-1)^n\ \nabla_F^{[n+1]}(\id\wedge\sigma_F)(\phi\tens\id) \cr
=&\ (\id\wedge\phi)(\id\wedge\sigma_E)(\nabla_E\tens\id)
+(\id\wedge\phi)\sigma_E(\id\tens\extd)  -\ (-1)^n\ \nabla_F^{[n+1]}(\id\wedge\sigma_F)(\phi\tens\id)   \cr
=&\ (\id\wedge(\id\wedge\phi)\sigma_E)(\nabla_E\tens\id)
+(\id\wedge\phi)\sigma_E(\id\tens\extd)  -\ (-1)^n\ \nabla_F^{[n+1]}(\id\wedge\sigma_F)(\phi\tens\id)   .
\end{align*}
Taking some care over the signs and using (\ref{severalsign}),
\begin{align*}
&\kern-20pt(-1)^n\ (\id\wedge\doublenabla(\phi))\sigma_E 
= (\id\wedge(\id\wedge\sigma_F)(\phi\tens\id))(\nabla_E\tens\id) \cr
&\ +\ (-1)^n\ (\id\wedge\sigma_F)(\phi\tens\id)(\id\tens\extd) 
 - \nabla_F^{[n+1]}(\id\wedge\sigma_F)(\phi\tens\id)   \cr
=&\ (\id\wedge\sigma_F)(\id\wedge\phi\tens\id)(\nabla_E\tens\id)
+(-1)^n\ (\id\wedge\sigma_F)(\phi\tens\extd) \cr
&\ -\ (\extd\wedge\sigma_F)(\phi\tens\id) 
-(-1)^n\  (\id\wedge\nabla_F^{[1]}\sigma_F)(\phi\tens\id)   \cr
=&\ (\id\wedge\sigma_F)(\id\wedge\phi\tens\id)(\nabla_E\tens\id)
+(-1)^n\ (\id\wedge\sigma_F)(\phi\tens\extd) \cr
&\ -\ (\id\wedge\sigma_F)((\extd\tens\id)\phi\tens\id)  - (-1)^n\  (\id\wedge(\id\wedge\sigma_F)(\nabla_F\tens\id))(\phi\tens\id)   \cr
&\ -\ (-1)^n\  (\id\wedge\sigma_F(\id\tens\extd))(\phi\tens\id) \cr
=&\ (\id\wedge\sigma_F)\big(      (\id\wedge\phi\tens\id)(\nabla_E\tens\id)    -\ ((\extd\tens\id)\phi\tens\id) 
-(-1)^n\  (\id\wedge\nabla_F\tens\id))(\phi\tens\id)   \big) \cr
=&\ -\ (\id\wedge\sigma_F)(\doublenabla(\phi)\tens\id)\ . \qquad\square
\end{align*}

\goodbreak
\subsection{Monoidal structure on ${}_A\mathcal{G}_A$}\label{secAGAmon}

\begin{theorem} \label{tensaga}
There is an associative functor $\boxtimes$ from ${}_A\mathcal{G}_A\times {}_A\mathcal{G}_A$ to ${}_A\mathcal{G}_A$, with product of objects being the usual tensor product of bimodules with bimodule connection. 
The tensor product of morphisms $\phi\in \mathrm{Mor}_n(E,G)$ 
and $\psi\in  \mathrm{Mor}_m(F,H)$ in ${}_A\mathcal{G}_A$ is
\begin{eqnarray*}
\phi\boxtimes\psi = (\id\wedge\sigma_F\tens\id)(\phi\tens\psi)
: E\otimes_A F \to\Omega^{n+m}\otimes_A G\otimes_A H.
\end{eqnarray*}
There is a signed rule for composition of tensor products:
\begin{eqnarray*}
(\phi\boxtimes\kappa)\circ(\psi\boxtimes\tau)\,=\,(-1)^{|\psi|\,|\kappa|}\
(\phi\circ\psi)\boxtimes(\kappa\circ\tau).
\end{eqnarray*}
For the differential of morphisms,
\begin{eqnarray*}
\doublenabla(\phi\boxtimes\psi) \,=\,\doublenabla(\phi)\boxtimes\psi +(-1)^n\,\phi\boxtimes\doublenabla(\psi) \ ,
\end{eqnarray*}
and for curvatures $R_{E\tens F} = R_{E}\boxtimes \id + \id\boxtimes R_F$. 
\end{theorem}
\proof First show that
given objects $(E,\nabla_E,\sigma_E)$ and $(F,\nabla_F,\sigma_F)$ then
$(E\tens_A F,\nabla_{E\tens F} ,\sigma_{E\tens F} )$ is also an object.
The extendable connection condition for $E\tens_A F$ follows fairly immediately from rearranging the crossings for $E\tens_A F$ (on the LHS of the following picture) in terms of the crossings for $E$ and the crossings for $F$ on the RHS:

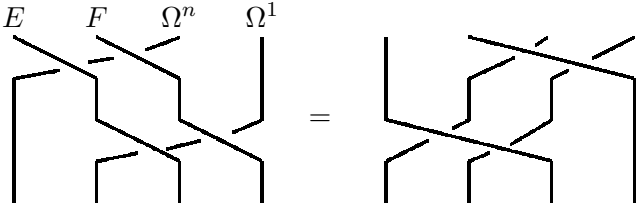
\unitlength 0.55 mm
\begin{picture}(160,50)(0,40)
\linethickness{0.3mm}
\multiput(10,80)(0.24,-0.12){83}{\line(1,0){0.24}}
\linethickness{0.3mm}
\multiput(30,80)(0.24,-0.12){83}{\line(1,0){0.24}}
\linethickness{0.3mm}
\put(70,70){\line(0,1){10}}
\linethickness{0.3mm}
\put(10,60){\line(0,1){10}}
\linethickness{0.3mm}
\put(30,60){\line(0,1){10}}
\linethickness{0.3mm}
\put(50,60){\line(0,1){10}}
\linethickness{0.3mm}
\put(70,60){\line(0,1){10}}
\linethickness{0.3mm}
\multiput(50,60)(0.24,-0.12){83}{\line(1,0){0.24}}
\linethickness{0.3mm}
\multiput(30,60)(0.24,-0.12){83}{\line(1,0){0.24}}
\linethickness{0.3mm}
\put(10,50){\line(0,1){10}}
\linethickness{0.3mm}
\put(100,60){\line(0,1){20}}
\linethickness{0.3mm}
\multiput(120,80)(0.48,-0.12){83}{\line(1,0){0.48}}
\linethickness{0.3mm}
\put(120,60){\line(0,1){10}}
\linethickness{0.3mm}
\put(140,60){\line(0,1){10}}
\linethickness{0.3mm}
\put(160,60){\line(0,1){10}}
\linethickness{0.3mm}
\multiput(100,60)(0.48,-0.12){83}{\line(1,0){0.48}}
\linethickness{0.3mm}
\put(160,50){\line(0,1){10}}
\linethickness{0.3mm}
\put(140,40){\line(0,1){10}}
\linethickness{0.3mm}
\put(160,40){\line(0,1){10}}
\linethickness{0.3mm}
\put(120,40){\line(0,1){10}}
\linethickness{0.3mm}
\put(100,40){\line(0,1){10}}
\linethickness{0.3mm}
\put(70,40){\line(0,1){10}}
\linethickness{0.3mm}
\put(50,40){\line(0,1){10}}
\linethickness{0.3mm}
\put(30,40){\line(0,1){10}}
\linethickness{0.3mm}
\put(10,40){\line(0,1){10}}
\linethickness{0.3mm}
\multiput(42,77)(0.32,0.12){25}{\line(1,0){0.32}}
\linethickness{0.3mm}
\multiput(28,74)(0.75,0.12){8}{\line(1,0){0.75}}
\linethickness{0.3mm}
\multiput(10,70)(0.65,0.12){17}{\line(1,0){0.65}}
\linethickness{0.3mm}
\multiput(63,57)(0.28,0.12){25}{\line(1,0){0.28}}
\linethickness{0.3mm}
\multiput(30,50)(0.59,0.12){17}{\line(1,0){0.59}}
\linethickness{0.3mm}
\multiput(47,53)(0.47,0.12){17}{\line(1,0){0.47}}
\linethickness{0.3mm}
\multiput(150,75)(0.24,0.12){42}{\line(1,0){0.24}}
\linethickness{0.3mm}
\multiput(120,70)(0.26,0.12){42}{\line(1,0){0.26}}
\linethickness{0.3mm}
\multiput(136,78)(0.18,0.12){17}{\line(1,0){0.18}}
\linethickness{0.3mm}
\multiput(140,70)(0.24,0.12){17}{\line(1,0){0.24}}
\linethickness{0.3mm}
\multiput(130,54)(0.2,0.12){50}{\line(1,0){0.2}}
\linethickness{0.3mm}
\multiput(120,50)(0.29,0.12){17}{\line(1,0){0.29}}
\linethickness{0.3mm}
\multiput(116,58)(0.24,0.12){17}{\line(1,0){0.24}}
\linethickness{0.3mm}
\multiput(100,50)(0.24,0.12){42}{\line(1,0){0.24}}
\put(10,85){\makebox(0,0)[cc]{$E$}}

\put(30,85){\makebox(0,0)[cc]{$F$}}

\put(50,85){\makebox(0,0)[cc]{$\Omega^n$}}

\put(70,85.5){\makebox(0,0)[cc]{$\Omega^1$}}

\put(84,60){\makebox(0,0)[cc]{$=$}}
\end{picture}
\medskip

Next we check the displayed condition in Lemma~\ref{mor2help} for $R_{E\tens F}$ to be a bimodule map (with   $n=1$ there):
 \begin{align*}
&\kern-20pt\nabla_{E\tens F}^{[1]}\sigma_{E\tens F} = (\extd\tens\id\tens\id)\sigma_{E\tens F}
-(\id\wedge\nabla_{E\tens F})\sigma_{E\tens F} \cr
=&\  (\extd\tens\id\tens\id)\sigma_{E\tens F}
-(\id\wedge\nabla_E\tens\id)\sigma_{E\tens F}    -\ (\id\wedge\sigma_E\tens\id)(\id\tens\id\tens\nabla_F)\sigma_{E\tens F} \cr
=&\  (\id\wedge\sigma_E\tens\id)(\nabla_E\tens\id\tens\id)(\id\tens\sigma_F)  + (\sigma_E\tens\id)(\id\tens\extd\tens\id)(\id\tens\sigma_F)    \cr
&\  -\ ((\id\wedge\sigma_E)(\sigma_E\tens\id)\tens\id)(\id\tens\id\tens\nabla_F)
(\id\tens\sigma_{F}) \cr
=&\  (\id\wedge\sigma_{E\tens F})(\nabla_E\tens\id\tens\id)  +  (\sigma_E\tens\id)(\id\tens\extd\tens\id)(\id\tens\sigma_F)    \cr
&\  -\ (\sigma_E\tens\id)(\id\tens\id\wedge\nabla_F)
(\id\tens\sigma_{F}) \cr
=&\  (\id\wedge\sigma_{E\tens F})(\nabla_E\tens\id\tens\id) + (\sigma_E\tens\id)(\id\tens\nabla_F^{[1]}\sigma_F)    \cr
=&\  (\id\wedge\sigma_{E\tens F})(\nabla_E\tens\id\tens\id) + (\sigma_E\tens\id)(\id\tens(\id\wedge\sigma_F)(\nabla_F\tens\id))    \cr
&\ +\  (\sigma_E\tens\id)(\id\tens \sigma_F(\id\tens\extd))    \cr
=&\  (\id\wedge\sigma_{E\tens F})(\nabla_{E\tens F}\tens\id) +\sigma_{E\tens F}(\id\tens\id\tens\extd).
\end{align*}

For morphisms the formula for $\phi\boxtimes\psi$ shows that it is a bimodule map, and we also check that it satisfies the $\sigma$ condition to be a morphism, which is best proved diagramatically:

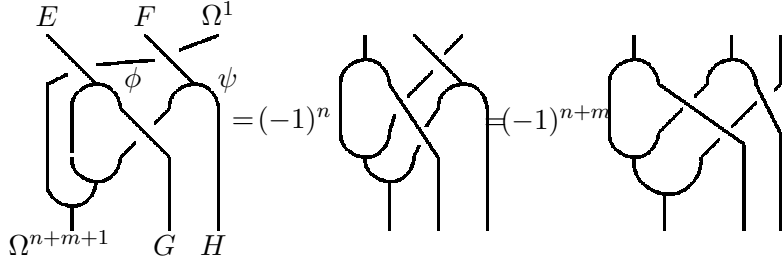
\unitlength 0.65 mm


\noindent
At each stage we have used the morphism property for one of $\phi$ and $\psi$, and also associativity of $\wedge$ product. The diagram for the property used for a single morphism appears at the end of the proof of Lemma~\ref{mor2help}, in the special case of the curvature (as $(-1)^2=1$ the sign does not appear there). 

 We also have to check associativity of $\boxtimes$ and the signed rule for compatibility with composition, both of which are routine exercises using the definitions and associativity of the wedge product (and best done diagrammatically). The more complicated case is the composition of tensor products, which is represented by the following diagram, where we use (\ref{severalsign}) and extendability:
 
 \unitlength 0.6 mm


Now we verify the formula for the derivative of a tensor product of morphisms: Suppose
${}_A\mathcal{G}_A$, for $\phi\in \mathrm{Mor}_n(E,G)$
and $\psi\in  \mathrm{Mor}_m(F,H)$,
Then
\begin{align*}
&\kern-20pt\nabla_{G\tens H}{}^{[n+m]}\circ(\phi\boxtimes\psi) =  (\extd\wedge\sigma_G\tens\id)(\phi\tens\psi)
+(-1)^n\ (\id\wedge\nabla_G{}^{[m]}\sigma_G\tens\id)(\phi\tens\psi) \cr
&\  +\ (-1)^{n+m}\ (\id\wedge\sigma_G\tens\id)(\phi\tens(\id\wedge\nabla_H)\psi) \cr
=&\  (\extd\wedge\sigma_G\tens\id)(\phi\tens\psi)
+(-1)^n\ (\id\wedge(\id\wedge\sigma_G)(\nabla_G\tens\id)\tens\id)(\phi\tens\psi) \cr
&\  +\ (-1)^n\ (\id\wedge\sigma_G(\id\tens\extd)\tens\id)(\phi\tens\psi)\cr
&\  +\ (-1)^{n+m}\ (\id\wedge\sigma_G\tens\id)(\phi\tens(\id\wedge\nabla_H)\psi) \cr
=&\  (\id\wedge\sigma_G\tens\id)(\nabla_G{}^{[n]}\phi\tens\psi
+(-1)^n\ \phi\tens\nabla_H{}^{[m]}\psi),
\end{align*}
while
\begin{align*}
&\kern-10pt\big(\id\wedge(\phi\boxtimes\psi)\big)\nabla_{E\tens F} =  
\big(\id\wedge(\phi\boxtimes\psi)\big)(\nabla_{E}\tens\id) +
\big(\id\wedge(\phi\boxtimes\psi)\big)(\sigma_E\tens\id)(\id\tens\nabla_{F})  \cr
=&\  (\id\wedge\sigma_G\tens\id)((\id\wedge\phi)\nabla_E\tens\psi) + \big(\id\wedge(\id\wedge\sigma_G\tens\id)(\phi\tens\psi)\big)
(\sigma_E\tens\id)(\id\tens\nabla_{F})  \cr
=&\  (\id\wedge\sigma_G\tens\id)((\id\wedge\phi)\nabla_E\tens\psi) + (-1)^n\ (\id\wedge\sigma_G\tens\id)(\phi\tens(\id\wedge\psi)\nabla_F).
\end{align*}
Subtracting these gives $\doublenabla(\phi\boxtimes\psi)$.  
\qquad$\square$

\medskip
We give four applications of this machinery. The first application is find another of the symmetries of the Riemann curvature,  this time involving a metric, written as $g\in\Omega^1 \tens_A\Omega^1$.

\begin{corollary} \label{yyuupp}
Suppose that $(\Omega^1,\nabla_{\Omega^1},\sigma_{\Omega^1})\in {}_A\mathcal{G}_A$ and that 
$g\in\Omega^1 \tens_A\Omega^1$ is preserved in the sense 
$\nabla_{\Omega^1\tens\Omega^1}g=0$. Then
\[
(R_{\Omega^1}\tens\id+(\sigma_{\Omega^1}\tens\id)(\id\tens R_{\Omega^1}))g=0.
\]
(Riemann antisymmetry identity.)
\end{corollary} 
\proof Here $\nabla_{\Omega^1\tens\Omega^1}g=0$ is equivalent to $\doublenabla(g)=0$ when $g$ is viewed as a morphism in ${}_A\CG$. This in turn implies $\doublenabla(\doublenabla(g))=0$ which is equivalent to   $R_{\Omega^1\tens\Omega^1}g=0$ by Theorem~\ref{the-bcat}.  Since
 $(\Omega^1,\nabla_{\Omega^1},\sigma_{\Omega^1})\in {}_A\mathcal{G}_A$ we can now use 
 Theorem~\ref{tensaga} to write $R_{\Omega^1\tens \Omega^1} = R_{\Omega^1}\boxtimes \id + \id\boxtimes R_{\Omega^1}$
 which gives the answer. \endproof
 
To find the classical identity corresponding to this we use 
 \[
R_{\Omega^1}(\extd x^i) = -\tfrac12\,R^i{}_{knm}\,\extd x^n\wedge\extd x^m\tens \extd x^k.
\]
and then
\[
R_{\Omega^1\tens \Omega^1}(g_{ij}\,\extd x^i\tens\extd x^j)=
 -\tfrac12\,R_{iknm}\,\extd x^n\wedge\extd x^m\tens \extd x^k\tens\extd x^i
  -\tfrac12\,R_{iknm}\,\extd x^n\wedge\extd x^m\tens \extd x^i\tens \extd x^k
\]
so we get the classical symmetry  $R_{abcd}=-R_{bacd}$ of the Riemann curvature of the Levi-Civita connection. 

For our second application, 
 it is explained in  \cite{bbsheaf} how to take the cohomology of a left module with connection $(F,\nabla_F)$ with zero curvature, and that this has many features in common with classical sheaf cohomology. Here we shall introduce a version of noncommutative cup product which is a direct analogue of the cup product in classical cohomology. 

\begin{proposition} \label{cohommult}
Suppose that $(E,\nabla_E,\sigma_E)$ is an extendable left bimodule connection with zero curvature, and $(F,\nabla_F)$ is a left module connection with zero curvature. Then the map
\[
\id\wedge \sigma_E \tens\id:
\Omega^m\otimes_A E \otimes \Omega^n\otimes_A F \to \Omega^{n+m}\otimes_A E\otimes_A F
\]
induces a product  $\cup:\coH^m(A,E,\nabla_E)\tens \coH^n(A,F,\nabla_F) \to \coH^{n+m}(A,E\tens_A F,\nabla_{E\tens F})$. 
\end{proposition}
\proof We use Lemma~\ref{mor2help},
\begin{align*}
\nabla_E{}^{[n]}\sigma_E = (\id\wedge\sigma_E)(\nabla_E\tens\id)+\sigma_E(\id\tens\extd) : E\otimes_A\Omega^n
\to \Omega^{n+1}\otimes_A E
\end{align*}
to show that
\begin{align*}
\nabla_{E\tens F}{}^{[n+m]}(\id\wedge \sigma_E \tens\id)=(\id\wedge \sigma_E \tens\id)(\nabla_E{}^{[m]}\tens\id+(-1)^m\id\tens 
\nabla_F{}^{[n]}),
\end{align*}
and the result then follows from routine algebra. \eproof

Classically the cup product is often used in conjunction with a product of sheaves, which in the language here would be a left module map
$:E\tens_A F \to G$ intertwining the connections, to give a product 
\begin{equation}
\cup:\coH^m(A,E,\nabla_E)\tens \coH^n(A,F,\nabla_F) \to \coH^{n+m}(A,G,\nabla_{G}).
\end{equation}
given by composing the map from Proposition~\ref{cohommult} with the cohomology map given by the morphism. We now turn to a third   application, which utilises the centre $Z_A(E)$ of a bimodule $E$, whose elements commute with every $a\in A$.

\begin{proposition}\label{promisedphi} Let $A$ be an algebra with exterior algebra $\Omega$ and  $(L,\nabla_L,\sigma_L)$  
an extendable bimodule connection on a line module. There is an algebra map $\hat\Phi_L:Z_A(\Omega)\to Z_A(\Omega)$ defined by $\hat\Phi_L(\xi)\tens e=\sigma_L(e\tens\xi)$ and if $R_L$ is a bimodule map (so $(L,\nabla_L,\sigma_L)\in {}_A\CG_A$) then there is a unique $\omega_L\in Z_A(\Omega^2)$ such that
$R_L(e)=\omega_L\tens e$. If $(M,\nabla_M,\sigma_M)$ is another such line module then so is their tensor product and
$\omega_{L\tens_A M}$ is given by $\omega_L+\hat\Phi_L(\omega_M)$. \end{proposition}
\proof Here $\hat\Phi_L$ extends the (inverse of) the Fr\"olich map to higher degree. If $\xi\in Z_A(\Omega^n)$ then the map $e\mapsto \sigma_L(e\tens\xi)$ is a right module map from
$L$ to $ \Omega^n\tens_A L$, so by Lemma~\ref{inventelement} it is given by $e\mapsto \hat\Phi_L(\xi)\tens e$ for some $\hat\Phi_L(\xi)$. Uniqueness also tells us that $\hat\Phi_{L\tens_A M}=\hat\Phi_L\circ \hat\Phi_M$.  Now suppose that $R_L$ is a bimodule map. The existence and uniqueness of $\omega_L\in\Omega^2$ follows from applying Lemma~\ref{inventelement} in the right $A$-module map case to $R_L$. Taking this form of $R_L$ and knowing that it is a left module map tells us that $\omega_L$ is central. For the curvature of the tensor product  we use Theorem~\ref{tensaga} and our observation about $\hat\Phi_{L\tens M}$. \eproof
 
Recall that in classical electromagnetism the curvature is often simply quoted as a 2-form rather than an operator. The above generalise this to line modules with connection in ${}_A\CG_A$, giving the curvature as a 2-form $\omega_L$. A fourth application  is to take traces of operators on objects $E\in {}_A\CG_A$ equipped with a quantum metric. We are interested particularly in the `trace' of powers of the curvature. 

\begin{corollary}\label{curvtrace} Let $E\in {}_A\CG_A$ and suppose that $g:A\to E\tens_A E$  and its inverse $(\ ,\ ):E\tens_A E\to A$ are morphisms in ${}_A\CG_A$ and  that $\doublenabla(g)=0$. Then 
we have an element of de Rham cohomology 
\[
[(\id\tens(,))(R_E{}^n\tens\id)g(1)]\in \coH_\mathrm{dR}^{2n}(A).
\] 
\end{corollary}
\proof In fact we only need to know that $g$ is a bimodule map and $\doublenabla(g)=0$ to know that it is a morphism. Here $\doublenabla(g)=0$ is equivalent to $\nabla_{E\tens E}(g)=0$ for the corresponding $g\in E\tens_A E$. From the equation $((\ ,\ )\tens\id)(\id\tens g)=\id:E\to E$ we deduce
\[
(\doublenabla((,))\boxtimes\id)\circ(\id\boxtimes g)+((\ ,\ )\boxtimes\id)\circ(\id \boxtimes \doublenabla(g))=
\doublenabla(\id)=0
\]
so  $\doublenabla((,))=0$ also holds. Now by Theorem~\ref{the-bcat} and Theorem~\ref{tensaga}, for all $n\ge 0$
\[
\doublenabla((\id\tens(,))(R_E{}^n\tens\id)g)=0
\]
which implies the result. \endproof

In the $n=0$ case we call this trace the `metric dimension' $\underline{\rm dim}_E=(\ ,\ )(g)$ of $E$ since classically the metric dimension will just recover the fibre dimension. For example, if $\Omega^1\in{}_A\CG_A$ and $\nabla_{\Omega^1}$ is metric compatible as in Corollary~\ref{yyuupp} then the above implies that the quantum dimension is in the kernel of $\extd$ and hence constant if the differential calculus is connected. In the case of the line module with connection in Proposition~\ref{promisedphi} preserving a metric on $L$ as just discussed we get the class of $\omega_L{}^n$ times the metric dimension in $\coH_\mathrm{dR}^{2n}(A)$.  It is not clear, however, what these noncommutative cohomology classes  are measuring in general. A similar problem of interpretation applies in the case of twisted cycles, in which (\ref{ncycledef}) is modified by a twisting automorphism of the algebra. This twisting theory is based on 
 \cite{KuMuTu}, and was further investigated in \cite{HadKra,NeTu,CarNesNesRen,wagsph}.  Following the classical theory, at first sight it might be thought that the twisted theory would also simply give information on the $K$-theory class of the module, but this is not correct: it gives a pairing with an equivariant $K$-theory class which incorporates the twisting automorphism. A possible lesson is that once we start to alter the definition of trace or $n$-cycle we move away from the classical idea of characteristic class, and that further study would be required, for example, to find out just what the metric de Rham classes above really are classifying. 
 
 Although not our main topic, one can also define a geometric Dirac operator given a connection $\nabla_S$ on a `spinor' bundle $\CS$ and additional data which we called in \cite{BMdirac} a `Clfford action' $\la:\Omega^1\tens_A\CS\to \CS$ subject to some axioms in order to get close to Connes' concept of a spectral triple at an algebraic level. The Dirac operator here is just $\dirac=\la\circ\nabla_\mathcal{S}$. We are now more concerned with the geometry and just ask what properties need for $\la$.  Classically, of course, this is part of an action of the Clifford algebra with an associative product being represented in $\CS$. We do not necessarily suppose that but we keep as motivation the  standard representation of such a Clifford algebra by wedge and interior product on the exterior algebra $\Omega$. 

\begin{lemma} \label{dirsqu} Suppose that the Clifford action extends to a bimodule map
$\la:\Omega^2\tens_A \CS\to \CS$ such that
\[
\varphi(\xi\la(\eta\la s)) =\kappa\, (\xi,\eta)\,s+(\xi\wedge\eta)\la s
\]
for all $\xi,\eta\in \Omega^1$, for some fixed invertible left $A$-bimodule map $\varphi:\CS\to\CS$
and constant $\kappa$. Let $(\Omega^1,\nabla_{\Omega^1},\sigma_{\Omega^1})$ be  a left bimodule connection
and suppose that the torsion $T$ is a bimodule map. Then
\[\dirac^2=\, \la\circ \doublenabla(\la)\circ\nabla_\mathcal{S}+\varphi^{-1}\circ(\kappa(\ ,\ )\nabla_{\Omega^1\tens \mathcal{S}}\nabla_\mathcal{S} + \la\circ R_\mathcal{S}
+(T\la\id)\nabla_\mathcal{S} ) \]
\end{lemma}
\proof We suppose that  $(\Omega^1,\nabla_{\Omega^1},\sigma_{\Omega^1})$ is a left bimodule connection and note that 
\[
\doublenabla(\la)=\nabla_\mathcal{S}\circ\la-(\id\tens\la)\nabla_{\Omega^1\tens \mathcal{S}}: \Omega^1\otimes_A  \mathcal{S}\to\Omega^1\otimes_A \mathcal{S}
\]
\[
\dirac^2=\la\circ\nabla_\mathcal{S}\circ\la\circ\nabla_\mathcal{S}=\la\circ \doublenabla(\la)\circ\nabla_\mathcal{S} +
\la\circ (\id\tens\la)\nabla_{\Omega^1\tens \mathcal{S}}\nabla_\mathcal{S} 
\]
We also have under our assumptions that 
\begin{align*}
\varphi\circ\la\circ (\id\tens\la)\nabla_{\Omega^1\tens \mathcal{S}}\nabla_\mathcal{S}=&\ \kappa\,((\,,)\sigma_{\Omega^1}\tens\id)(\id\tens\nabla_\mathcal{S})\nabla_\mathcal{S}+\kappa\,((\,,)\nabla_{\Omega^1}\tens\id)\nabla_\mathcal{S}
\\
&+(\wedge\sigma_{\Omega^1}\tens\id)(\id\tens\nabla_\mathcal{S})\nabla_\mathcal{S}
+ (\wedge\nabla_{\Omega^1}\tens\id)\nabla_\mathcal{S}\end{align*}
where we expanded $\nabla_{\Omega^1\tens \mathcal{S}}$ in order particularly to apply $\wedge$. If we suppose that $\nabla_{\Omega^1}$ has torsion $T$ a bimodule map, so that $\wedge\sigma_{\Omega^1}=-\wedge$, then this simplifies to give
\begin{align*}
\varphi\circ\dirac^2=&\ \varphi\circ\la\circ \doublenabla(\la)\circ\nabla_\mathcal{S} +
\kappa\,((\,,)\sigma_{\Omega^1}\tens\id)(\id\tens\nabla_\mathcal{S})\nabla_\mathcal{S}
+\kappa\,((\,,)\nabla_{\Omega^1}\tens\id)\nabla_\mathcal{S}\\ & + \la\circ R_\mathcal{S}
 +(T\la\id)\nabla_\mathcal{S} 
\end{align*}
which gives the formula stated on recognising the torsion and curvature. \endproof

The first term in this result vanishes if  the Clifford action is covariant for the relevant connections so that $\doublenabla(\la)=0$, which is the case classically since it is defined by the metric and the relevant connections are compatible. We also normally ask for the torsion to vanish as we use the quantum Levi-Civita connection, leaving us a Laplacian-Beltrami operator $(\ ,\ )\nabla_{\Omega^1\tens \mathcal{S}}\nabla_\mathcal{S}$ and curvature $\la\circ R_\mathcal{S}$, much as in the classical Lichnerowicz formula. This Laplace-Beltrami operator on spinors  clearly generalises the geometric Laplacian naturally defined on $A$ itself by $\Delta=(\ ,\ )\nabla_{\Omega^1}\extd$, while $\la\circ R_\mathcal{S}$ classically would be a multiple of the Ricci scaler and suggests a new `Clifford action' approach to the contraction needed for the Ricci tensor and Ricci scaler. Such further developments will be considered elsewhere, particularly in but not limited to the category ${}_A\CG_A$.

\section{Examples} \label{secex}

We collect some examples at the end of the paper in order not to distract from the general results. Since our results are very general there are plenty more examples, although for the Chern theory we have the same limitations as in cyclic cohomology due to the lack of cycles of the relevant dimension (eg for the $q$-sphere we would need twisted cycles as in \cite{KuMuTu} which requires more analysis for our results). We also discover that in some noncommutative models, extendability is linked to zero curvature. 

\subsection{Connections on $\Omega^1(S_3)$}

Finite groups $G$ provide nice examples where the differentials are noncommutative even though coordinate algebra $A=\kk(G)$ is commutative.
A bicovariant calculus here is fixed by $\CC\subseteq G\setminus\{e\}$ stable under conjugation and has left invariant basis elements $e_a$ for $a\in \CC$ and relations and exterior derivative 
\[ e_a.f=R_a(f)e_a,\quad \extd
f=\sum_{a\in {\mathcal{C}}}(R_a(f)-f) e_a\]
where $f\in \kk (G)$ and $R_a(f)(x)=f(xa)$. The canonical `Woronowicz construction'\cite{Wor} exterior algebra in this case is defined via the crossed module braiding $\Psi(e_a\tens e_b)=e_{aba^{-1}}\tens e_{a}$ and in degree 2 consists in setting $\ker(\id-\Psi)$ to zero. The above is immediate from a general theory of Hopf algebras in \cite{Wor}, while a recent nontrivial result \cite{MaRie} is that the noncommutative de Rham cohomology for the calculus is $H^1(G)=\kk\theta$ for $\kk$ with $\theta=\sum_ae_a$, when $\kk$ has characteristic zero and $\CC$ as a quandle is `locally skew' (a concept that includes all Weyl groups of complex semisimple Lie algebras with $\CC$ given by reflections). We refer to \cite{MaRie} for details. We use this setting to illustrate our above results on $G=S_3$, the permutation group of 3 objects, and $\CC=\{u,v,w\}$ the set of 2-cycles (or transpositions) in $S_3$ where $u^2=v^2=e$ and $w=uvu=vuv$. The quantum metric here is the Euclidean one $g=e_u\tens e_u+e_v\tens e_v+e_w\tens e_w$.

(i) Then the canonical exterior algebra $\Omega(S_3)$ in for this calculus  is well-known to have relations
\[ e_u\wedge e_v+e_v\wedge e_w+e_w\wedge e_u=0,\quad e_v\wedge e_u+e_u\wedge e_w+e_w\wedge e_v=0,\quad e_u^2=e_v^2=e_w^2=0\]
\[ \extd e_u+e_v\wedge e_w+e_w\wedge e_v=0,\quad \extd e_v+e_w\wedge e_u+e_u\wedge e_w=0,\quad \extd e_w+e_u\wedge e_v+e_v\wedge e_u=0\]
which imply dimensions $1:3:4:3:1$ and in particular a unique up to scale `top form'
\[ {\rm Vol}=e_u\wedge e_v\wedge e_u\wedge e_w=e_v\wedge e_u\wedge e_v\wedge e_w=-e_w\wedge e_u\wedge e_v\wedge e_u=-e_w\wedge e_v\wedge e_u\wedge e_v   \]
and equal to the two cyclic permutations $u\to v\to w\to u$. 

We observe that this exterior algebra has a well-defined $4$-cycle  defined by $\int(f{\rm Vol})=\int f=\sum_{g\in S_3}f(g)$. First, from the above versions of ${\rm Vol}$ we see that swapping $e_u\wedge e_v\wedge e_u$ with $e_w$ gives  a minus sign and similarly for other identities for (\ref{ncycledef}) for $\omega,\rho$ 1- and 3-forms with constant coefficients. Here $e_u\wedge e_v\wedge e_w=e_w\wedge e_v\wedge e_u=-e_w\wedge e_u\wedge e_w=-e_u\wedge e_w\wedge e_u$ and its two cyclic permutations are a convenient basis of 3-forms. Similarly in degree 2 we have basis $e_u\wedge e_v, e_v\wedge e_u, e_v\wedge e_w, e_w\wedge e_v$ and these mutually commute so that (\ref{ncycledef}) holds on basic 2-forms as well. Also note that ${\rm Vol}$ is zero unless the total degree given by the product of the basis labels is $e$ as it is known that ${\rm Vol}$ here is central (these basic facts were used recently in \cite{Ma:bfou}). Then for the general case (where we sum over labels belonging to our declared bases) 
\[ \int(\omega\wedge \rho)=\int \omega_a e_a\wedge \rho_{bcd} e_b\wedge e_c\wedge e_d=\int \omega_a R_{a}( \rho_{bcd})e_a\wedge  e_b\wedge e_c\wedge e_d\]
\[=-\int R_{bcd}(\omega_a)  \rho_{bcd}e_b\wedge e_c\wedge e_d\wedge e_a  =-\int  \rho_{bcd}  e_b\wedge e_c\wedge e_d \wedge \omega_a e_a=-\int\rho\wedge\eta\]
for degrees 1,3 (and similarly for degrees 3,1), and
\[ \int(\omega\wedge \rho)=\int \omega_{ab} e_a\wedge e_b\wedge \rho_{cd}  e_c\wedge e_d=\int \omega_{ab} R_{ab}( \rho_{cd})e_a\wedge  e_b\wedge e_c\wedge e_d\]
\[=\int R_{cd}(\omega_{ab})  \rho_{cd} e_c\wedge e_d\wedge e_a\wedge e_b     =\int  \rho_{cd}   e_c\wedge e_d \wedge \omega_{ab} e_a\wedge e_b=\int\rho\wedge\eta\]
for degree 2,2. We used of course that $\sum$ over $S_3$ is translation invariant. We have $\int\extd \omega=\int[\theta,\omega\}=0$ (where we use the graded-commutator) since the calculus is inner by $\theta=e_u+e_v+e_w$.

(ii) Next we look at connections on $\Omega^1$. We focus on left-invariant connections,  which means of the form
\[
\nabla_{\Omega^1} e_a=-\sum_{b,c\in\mathcal{C}}\Gamma^a{}_{bc}\, e_b\tens e_c
\]
where $\Gamma^a{}_{bc}\in \kk$. For simplicity we further restrict to the case where $\nabla_{\Omega^1}$ is also right-invariant, which reduces to $\Gamma^{gag^{-1}}{}_{gbg^{-1},gcg^{-1}}=\Gamma^{a}{}_{b,c}$ for all $a,b,c\in \mathcal{C}$ and $g\in G$ in the case of general group. In our case, as  conjugation induces all possible permutations of the set $\mathcal{C}$,
we deduce that there are
only 5 possible different values of the Christoffel symbols $\Gamma^a{}_{bc}$, namely of the form 
$\Gamma^x{}_{xx}$, $\Gamma^x{}_{yz}$, $\Gamma^x{}_{yx}$,
$\Gamma^x{}_{xy}$ and $\Gamma^y{}_{xx}$ (where $x,y,z$ are all different). We set 
\begin{eqnarray} \label{chr5}
\Gamma^x{}_{xx}=a-1,\quad \Gamma^x{}_{yz}=c,\quad \Gamma^x{}_{yx}=d-1,\quad
\Gamma^x{}_{xy}=e,\quad \Gamma^y{}_{xx}=b\end{eqnarray}
with parameters  $a,b,c,d,e\in\kk
$   (these should not be confused with generic elements of $\mathcal{C}$ which we won't need at this point) and since we use computer algebra we assume $\kk$ has characteristic zero. The curvature in our context is
\begin{align*}
R_{\Omega^1}(e_a)=-&\ \big(\Gamma^a{}_{bc}\extd e_b+\Gamma^a{}_{bs}\Gamma^s{}_{rc}
e_b\wedge e_r\big)\tens e_c .
\end{align*}
which in terms of the parameters works out as 
The curvature is given by the following (and permutations....)
\begin{align*}
R_{\Omega^1}&(e_u) =e_w\wedge e_u \tens \big( (-b^2 - a d + d^2 + c e)e_u+(d - c) (b - e)e_v+(b^2 + a c - c^2 - d e)e_w
\big) \cr
&+e_u\wedge e_v \tens \big( (b c - a d + d^2 - e^2)e_u+(-a b + b d + c d - a e - b e + c e)e_v+(b - e) (b + d + e)e_w
\big) \cr
&+e_u\wedge e_w \tens
\big( (b c - a d + d^2 - e^2)e_u+(b - e) (b + d + e)e_v+(-a b + b d + c d - a e - b e + c e)e_w
\big) \cr
&+e_v\wedge e_u \tens
\big((-b^2 - a d +d^2 + c e)e_u+(b^2 +a c - c^2 - d e)e_v+(d - c) (b - e)e_w
\big) \ .
\end{align*}
Using the relations and monomial equalities in (i) one can check that the first Bianchi identity holds in the form \[ \wedge R_{\Omega^1}(e_u)=0\]
and similarly for $e_v,e_w$ as it must in the absence of torsion by Lemma~\ref{1stbia}. The 2nd Bianchi identity 
\[ (\extd\tens\id +\id\wedge\nabla_{\Omega^1}) R_{\Omega^1}=(\id\wedge R_{\Omega^1})\nabla_{\Omega^1}\]
in Theorem~\ref{the-bcat}  must also hold and can be verified.  A short compatutation gives 
\[ (\id\tens(\ ,\ ))(R_{\Omega^1}\tens\id)g={\rm Tr}R_{\Omega^1}=0\]
where the metric trace and usual trace coincide because of the trivial Euclidean form of the metric. The same is true for the trace of $R_{\Omega^1}{}^2$ but we check it in the following more manageable case.

(iii) Geometrically speaking, we are interested in ad-invariant connections that are torsion free, which is $c-1=e=d-1$, and cotorsion-free, which is $c-1=b=d-1$. These are as close as one can come to `Levi-Civita' in this example with respect to the Euclidean metric $g=\sum_a e_a\tens e_a$. There is a two-parameter moduli space of such connections cf\cite{Ma:non} and we write the parameters as $\lambda= c-a-2$ and $\mu=-c$. Then one can compute that
\begin{align*} R_{\Omega^1}(e_u)=&\ - \extd e_u\tens e_u-\extd e_v\tens e_w-\extd e_w\tens e_v\\
&-(3+\lambda)\big(\mu\,\extd e_u\tens e_u+(2(1+\mu)e_u\wedge e_v-\mu e_v \wedge e_u)\tens e_v \big)\\
&-(3+\lambda)\big(2(1+\mu)e_u \wedge e_w-\mu e_w \wedge e_u)\tens e_w\big)\end{align*}
with the same under $e_u\to e_v\to e_w\to e_u$. We write the curvature in the shorthand form 
\[ R_{\Omega^1}(e_u)=\alpha (e_v\wedge e_w+ e_w\wedge e_v)\tens e_u+(\beta e_u\wedge e_v+\gamma e_v\wedge e_u)\tens e_v+ (\beta e_u\wedge e_w+\gamma e_w\wedge e_u)\tens e_w)\]
for coefficients $\alpha=(3+\lambda)\mu+1, \beta=-1-2(3+\lambda)(1+\mu), \gamma=(3+\lambda)\mu-1$ and the same under $u\to v\to w\to u$. Then trace of $R_{\Omega^1}{}^2$ is given by the $\tens e_a$ coefficient of $R_{\Omega^1}{}^2(e_a)=(\wedge\tens\id)(\id\tens R_{\Omega^1})R_{\Omega^1}(e_a)$,  summed over $a=u,v,w$. In fact each of these terms vanish, 
\begin{align*} R_{\Omega^1}{}^2(e_u)|_{\tens e_u}=&\alpha^2(e_v\wedge e_w+ e_w\wedge e_v)^2+(\beta e_u\wedge e_v+\gamma e_v\wedge e_u)\wedge(\beta e_v\wedge e_u+\gamma e_u\wedge e_v)\\
&+(\beta e_u\wedge e_w+\gamma e_w\wedge e_u)\wedge(\beta e_w\wedge e_u+\gamma e_u\wedge e_w)
\end{align*}
where in the first term we apply $R_{\Omega^1}$ to the $\tens e_u$ term of $R_{\Omega^1}(e_u)$, for second term we apply it to the $\tens e_v$ term of $R_{\Omega^1}(e_u)$ and for third term we apply it to the  $\tens e_w$ terms of $R_{\Omega^1}(e_u)$,  and in all cases we pick off the $\tens e_u$ term of 
this second application of $R_{\Omega^1}$.  Now using the relations of the exterior algebra we obtain
\[ R_{\Omega^1}{}^2(e_u)|_{\tens e_u}=\alpha^2( (e_v \wedge e_w)^2  + (e_w\wedge e_v)^2)+\beta\gamma((e_u\wedge e_v)^2+(e_v\wedge e_u)^2+(e_u\wedge e_w)^2+(e_w\wedge e_u)^2)=0\]
since all of these squares vanish using the commutation relations (for example, $e_u\wedge e_v\wedge e_u\wedge e_v=-e_u\wedge(e_u\wedge e_w+e_w\wedge e_v)\wedge e_v=0$). Hence $\int \Tr R_{\Omega^1}{}^2=0$, at least for this class of connections, as each term in the trace vanishes. Here integration just sums over the group but all our coefficients are constant so this just multiplies by the order of the group.

(iv) We next ask when $(\Omega^1,\nabla_{\Omega^1})$ is an objects of ${}_A\CG_A$ for $A=C(S_3)$ with the same calculus.  From the bimodule commutation relations we have
\[
\big(\Gamma^a{}_{bc}\extd e_b+\Gamma^a{}_{bs}\Gamma^s_{rc}
e_b\wedge e_r\big)\, f=R_{ac^{-1}}(f)\, \big(\Gamma^a{}_{bc}\extd e_b+\Gamma^a{}_{bs}\Gamma^s{}_{rc}
e_b\wedge e_r\big)
\]
for any $f$ and the $\tens e_c$ component of $R_{\Omega^1}(e_a)$ for a connection in setting above. In view of this, a calculation for our full 5-parameter moduli space of bicovariant connections gives the following conditions on the parameters for the curvature being a right- and hence bi-module map:

(1)\quad $b=c=d=e=0$,

(2)\quad $a=-b(b+c)/c$, $c\ne 0$, $d=c$, $e=-b-c$, 

(3)\quad $d=a$,  $e= c=b$,

(4)\quad $a=c$,  $e=d=b$,

(5)\quad $a=b+c-b^2/c$, $c\ne 0$, $d=c$, $e=b$. 

\smallskip\noindent of which one can check using the formula for $R_\nabla$ in part (ii) that (1)-(4) are just the cases of zero curvature and the only possibly nonzero case of the curvature a bimodule map is (5),  where 
\[ R_\nabla(e_u)={ (c-e)^2(c+2e)\over c}(e_u\wedge e_v\tens e_v+ e_u\wedge e_w\wedge e_u).\]
Here $e=c,-c/2$ fall back to an instance of cases (4),(2) respectively. 

Moreover, it is  known \cite{BMriem} that all left connections on $\Omega^1$ as in (ii) are bimodule ones and using the tensor product basis $\{e_I\}$ in the lexicographical order $e_u\tens e_u$, $e_u\tens e_v$, $e_u\tens e_w$, $e_v\tens e_u$ to $e_w\tens e_w$ that
{\renewcommand{\arraystretch}{1.3}
\[\sigma= {\small\small \small \begin{pmatrix}
 a & 0 & 0 & 0 & b & 0 & 0 & 0 & b \\
 0 & e & 0 & 0 & 0 & c & d & 0 & 0 \\
 0 & 0 & e & d & 0 & 0 & 0 & c & 0 \\
 0 & 0 & c & e & 0 & 0 & 0 & d & 0 \\
 b & 0 & 0 & 0 & a & 0 & 0 & 0 & b \\
 0 & d & 0 & 0 & 0 & e & c & 0 & 0 \\
 0 & c & 0 & 0 & 0 & d & e & 0 & 0 \\
 0 & 0 & d & c & 0 & 0 & 0 & e & 0 \\
 b & 0 & 0 & 0 & b & 0 & 0 & 0 & a
\end{pmatrix}}
\]
}
where $\sigma_{\Omega^1}(e_I)=\sum_J  \sigma_I{}^Je_J$ and $\{\sigma_I{}^J\}$ is the matrix as shown (for example $\sigma_{\Omega^1}(e_u\tens e_v)=e\,e_u\tens e_v+c\,e_v\tens e_w+d\,e_w\tens e_u$ according to the second row).  

We also ask which of our 5-parameter moduli of connections is  extendable. Here the extended  $\sigma_{\Omega^1}:\Omega^1\tens_A\Omega^2\to\Omega^2\tens_A\Omega^1$, if it exists, must obey 
\[
\sigma_{\Omega^1}(\omega\tens \eta\wedge \zeta)=(\wedge\tens\id)(\id\tens\sigma_{\Omega^1})(\sigma_{\Omega^1}\tens\id)(\omega\tens \eta\tens \zeta)
\]
and if so then it further extends to higher forms to give an extendable bimodule connection in the sense of Definition~\ref{extendableconn} since the relations of the exterior algebra in the present example are quadratic. We ask when this is well defined. For example
\begin{align*} 0&=\sigma_{\Omega^1}(e_u\tens e_u\wedge e_u)=(\wedge\tens\id)(\id\tens\sigma_{\Omega^1})(a e_u\tens e_u+ b e_v\tens e_v+ b e_w\tens e_w)\tens e_u)\\
&=ae_u\wedge (be_v\tens e_v+ b e_w\tens e_w)+be_v\wedge(c e_u\tens e_w+d e_w\tens e_v)+b e_w\wedge(c e_u\tens e_v+d e_v\tens e_w)
\end{align*}
which given the quadratic relations of the algebra requires $ba=bd=bc$. Similarly $0=\sigma_{\Omega^1}(e_u\tens e_v\wedge e_v)$ gives us  $e a= e d=e c$ and $cd=e b$. Proceeding in the same way, $\sigma_{\Omega^1}(e_u\tens( \ ))$ applied to the two 3-term relations give us 
\[ e^2+ad+ec=b^2+ad+bc=d^2+bc+ce,\quad e^2+ac+ed=b^2+ac+bd=c^2+bd+de\]
Taken together, these nine equations have six 1-parameter moduli of solutions as follows:

(1) $b=c=d=e=0$   

(2) $a=c=d$, $b=q^{-1} a$, $e=q a$ where $q$ is any cube root of 1 (three cases)

(3) $a=d$, $b=c=e=0$

(4) $a=c$, $b=d=e=0$.

\noindent We have done this by hand but a general approach is to  ask which $\sigma_{\Omega^1}$ obey the mixed braid relation
\begin{align} \label{pu67}
(\id\tens\sigma_{\Omega^1})(\sigma_{\Omega^1}\tens\id)(\id\tens\Psi)=(\Psi\tens\id)(\id\tens\sigma_{\Omega^1})(\sigma_{\Omega^1}\tens\id),
\end{align}
where $\Psi(e_u\tens e_v)=e_w\tens e_u$ etc. is the crossed-module braiding that defines the calculus quadratic relations in our case. 
Such $\sigma_{\Omega^1}$ extend to the standard exterior algebra at least to degree 2 where the relations are given by $\ker(\id-\Psi)$. 
In the present example, we obtain all extendable $\sigma_{\Omega^1}$, i.e. the same moduli (1)-(4) as found directly. 
On the other hand, one can easily see that these cases are all instances of the corresponding moduli of $R_{\Omega^1}=0$ connections already computed: 

\begin{lemma} All six 1-parameter moduli of ad-invariant extendable connections make $\Omega^1$ a flat object of ${}_{C(S_3)}\mathcal{G}_{C(S_3)}$.
\end{lemma}

Exactly two points among all six moduli are cotorsion-free, namely
\[ a=c=d=b+1={1\over 1-q^{-1}},\quad e=q d\]
where $q\ne 1$ is a cube root of unity, but neither of these is torsion-free. Exactly two points among all six moduli are torsion free, namely 
\[ a=c=d=e+1={1\over 1-q},\quad b=q^{-1} d\]
where $q\ne 1$ is a cube root of unity, but neither of these is cotorsion-free.

(iv) We can also ask about full metric compatibility. For the full 5-parameter moduli of connections, the equations of $\nabla_{\Omega^1\tens\Omega^1}g=0$ of metric compatibility come down to
\[
 a^2+2 b e=1\ ,\ b c + e d + c d=0\ ,\  b e + c^2 + d^2=1\ ,\  a b + a e + b e=0\ ,\ e c + b d + c d=0.
\]
This has a discrete moduli of solutions as follows. 

\smallskip (a): $a=\pm1$, $b=d=e=0$, $c=\eps a$  where $\eps=\pm 1$.

\smallskip (b): $a=\pm 1$, $b=c=e=0$, $d=\eps a$ where $\eps=\pm 1$.

\smallskip (c): $a=c=d=\pm 1/\sqrt{3}$, $b=q^{-1}a$, $e=qa$ where $q\neq 1$ is a cube root of 1.

\smallskip (d): $a=\pm\, \mathrm{i}$, $b=za$, $c=d=0$, $e=-z^{-1}a$  where $z^2-z-1=0$  (for the golden ratio).

\smallskip (e): $a=\pm 1/3$, $b=e=-2a$, $c=xa$, $d={2xa \over x-2}$ where $(x-1) (x+2) (x^2 - 5 x +  10)=0$. 

\smallskip None of these are torsion free, so there are no strict quantum Levi-Civita connections as fore-warned already in part (ii).  

We also see that the metric preserving bimodule connections which are extendable  are the four cases (a),(b) both with $\eps=1$, and (c) with the two values of $q$. These are flat as we know from the above. The other curvatures are respectively,
\begin{align*} &(a)\quad \eps=-1:\quad R_{\Omega^1}(e_u)=-2(e_w\wedge e_u\tens e_w+e_v\wedge e_u\tens e_v)\\
& (b) \quad \eps=-1:\quad R_{\Omega^1}(e_u)=-2(e_v\wedge e_w+e_w\wedge e_u)\tens e_u\\
& (d)\quad R_{\Omega^1}(e_u)= z^2\left(e_w\wedge e_u\tens (e_u-e_w)+e_v\wedge e_u\tens (e_u-e_v)\right)\\
&\kern70pt+{1\over z^2}e_u\wedge \left(e_v\tens( e_u+(1-z^4)e_w)+e_w\tens e_u+(1-z^4)e_v\right) \\
& (e)\quad R_{\Omega^1}(e_u)= {(x-4) (x-1) (x+2)\over 9 (x-2)^2}\\
&\kern80pt \big(2(e_v\wedge e_w+ e_w\wedge e_v)\tens e_u+(2-x)(e_w\wedge e_u\tens e_w+e_v\wedge e_u\tens e_v)\big).\end{align*}

\subsection{Riemannian geometry of $\C_q[SL_2]$ with its 4D calculus}

It is well known that $\C_q[SL_2]$ has a natural 3D left covariant and 4D bicovariant calculus \cite{Wor} (the latter is unique if want a classical limit); the first does not admit a central quantum metric while the second has a unique invariant one 
\[g=e_c\tens e_b+q^2 e_b\tens e_c +{q^3\over(2)_q}(e_z\tens
e_z-\theta\tens\theta)\]
where we used the standard conventions as in \cite{Ma:bfou} with basic 1-forms $e_a,e_b,e_c,e_d$ (they arise as a $2\times 2$ matrix) with $q$-commutation relations with the coordinate algebra generators and $q$-anticommutation relations among themselves. We let $e_z=q^{-2}e_a-e_d$ so $e_b,e_c,e_d$ $q$-deform the classical 1-forms and there is an extra direction $\theta=e_a+e_b$ which has no classical meaning but which makes the calculus inner when $q\ne 1$. We let $\lambda=1-q^{-2}$. 

This quantum metric has a torsion-free and cotorsion-free or `weak quantum Levi-Civita' connection found  in \cite{Ma:ric} using a quantum frame bundle approach, namely
\begin{eqnarray*} \nabla_{\Omega^1} e_a&=&-\nabla_{\Omega^1} e_d
={1\over (2)_{q^2}}\left(e_b\tens
e_c-e_c\tens e_b- \lambda{q^3\over (2)_q}e_z\tens e_z\right)
\\
 \nabla_{\Omega^1} e_b&=& {1\over
(2)_{q^2}} \left(e_z\tens e_b-q^2 e_b\tens e_z\right),\quad \nabla_{\Omega^1} e_c=
{1\over (2)_{q^2}} \left(-q^2 e_z\tens e_c+
e_c\tens e_z\right)\end{eqnarray*}
and is unique at least for generic $q$. It is straightforward to check using the exterior algebra relations (in the form listed in \cite{Ma:bfou}) that its curvature
\begin{eqnarray*} R_{\Omega^1} e_a&=&-R_{\Omega^1}(e_d)=
-{q^2\over (2)_{q^2}^2}\left(e_c\wedge e_z\tens e_b+q^{-2}e_b\wedge e_z\tens e_c+\lambda
e_b\wedge e_c\tens e_z\right)\\
 R_{\Omega^1} e_b&=&{q^2\over (2)_{q^2}^2}\left(e_b\wedge e_z\tens
e_z+q^{-3}(2)_q e_b\wedge e_c\tens e_b\right)\\
R_{\Omega^1} e_c &=&{q^2\over (2)_{q^2}^2}\left(e_c\wedge e_z\tens
e_z-q^{-1}(2)_qe_b\wedge e_c\tens e_c\right).\end{eqnarray*}
has $ \wedge R_{\Omega^1}=0$ as required by the 1st  Bianchi identity. The 2nd Bianchi identity is more involved but must also hold. 

This $\nabla_{\Omega^1}$ is {\em not} a bimodule connection (one can determine what $\sigma_{\Omega^1}(e_a\tens e_b)$ etc would have to be and check that they are not compatible with the relations of the 4D calculus). So a `weak Levi-Civita' connection is the best we can do as we do not have a concept of full metric compatibility. This also means that the theory of Section~5 does not apply.  We can nevertheless still compute the quantum metric traces and they are
\[ \und\dim_{\Omega^1}=(2)_q^2,\quad  (\id\tens (\ ,\ ))(R_{\Omega^1}\tens\id)g=0,\quad (\id\tens (\ ,\ ))(R_{\Omega^1}{}^2\tens\id)g=0\]
which deforms or agrees with the classical picture in so far as one can compare. Note that the lack of a full quantum Levi-Civita connection here and only a weak torsion-free cotorsion-free one is similar to the $C(S_3)$ example above. They also have the same noncommutative de Rham cohomology and in both cases the quantum Ricci tensors are essentially proportional to the metric\cite{Ma:non,Ma:ric}. 

\subsection{Bicrossproduct model quantum spacetime revisited}\label{bicross}

Here we look at the 2D noncommutative model with curvature found in \cite{BegMa2} for the quantum spacetime $A$ with differential algebra 
\[ [r,t]=\lambda r,\quad [r,\extd t]=\lambda\extd r,\quad [t,\extd t]=\lambda\extd t,\quad [r,\extd r]=0=[t,\extd r]\]
with $\extd t,\extd r$ anticommuting, which admits a unique form of central quantum metric
\[ g=((1+b\lambda^2)\extd r-\lambda b v)\tens\extd r+b v\tens v,\quad  v=r\extd t- t\extd r.\]
Here  $\extd r,v$ are central and $\extd r\wedge v+v\wedge\extd r=0, (\extd r)^2=0$ but $v^2=\lambda {\rm Vol}$, where ${\rm Vol}=v\wedge\extd r$ is of top degree. It is already known \cite{BegMa2} that there are two moduli of $*$-preserving metric compatible bimodule connections for this model when $b$ is real and $\lambda$ imaginary, namely a line and a conic, with a unique torsion free `quantum Levi-Civita connection' in each component. Only one of these, on the conic, has a classical limit as $\lambda\to 0$ so there is a unique quantum Levi-Civita connection with classical limit. 

(i) Here we will look at bimodule connections more generally but still of the homogeneous form studied in \cite{BegMa2}, where
\[ \nabla_{\Omega^1}\extd r= {1\over r}(\alpha v\tens v+\beta v\tens\extd r+\gamma\extd r\tens v+\delta\extd r\tens\extd r)\]
\[ \nabla_{\Omega^1} v= {1\over r}(\alpha' v\tens v+\beta' v\tens\extd r+\gamma'\extd r\tens v+\delta'\extd r\tens\extd r)\]
for constants $\alpha,\cdots,\delta,\alpha',\cdots,\delta'$. These have braiding
$\sigma_{\Omega^1}(\extd r\tens \extd r)=\extd r\tens \extd r$, $\sigma_{\Omega^1}(v\tens \extd r)=\extd r\tens v$ and
\[ \sigma_{\Omega^1}(\extd r\tens v)=\lambda\alpha v\tens v+(1+\lambda\beta)v\tens\extd r+\lambda\gamma\extd r\tens v+\lambda\delta \extd r\tens\extd r   \]
\[  \sigma_{\Omega^1}(v\tens v)=(1+\lambda\alpha') v\tens v+\lambda\beta' v\tens\extd r+\lambda\gamma'\extd r\tens v+\lambda\delta' \extd r\tens\extd r,   \]
torsion \cite{BegMa2}
\[ T(\extd r)={1\over r}(\lambda\alpha+\beta-\gamma){\rm Vol},\quad T(v)={1\over r}(\lambda\alpha'+\beta'-\gamma'+2){\rm Vol})\]
and curvature of the form
\[ R_{\Omega^1}(\extd r)=-{1\over r^2}{\rm Vol}\tens (c_1 v+ c_2\extd r),\quad R_{\Omega^1}(v)=-{1\over r^2}{\rm Vol}\tens (c_3 v+ c_4\extd r)\]
for some constant coefficients $c_i$  given explicitly in \cite{BegMa2} in terms of our eight parameters. Both the torsion and curvature are bimodule maps. The Bianchi identities are also automatic as $\Omega^3=0$. The quantum dimension from\cite{BegMa2} and now the quantum metric trace of the curvature using
\[ (R_{\Omega^1}\tens\id)(g)=-{(1+b \lambda^2)\over r^2}{\rm Vol}\tens (c_1 v+c_2\extd r)\tens \extd r-{b\over r^2}{\rm Vol}\tens (c_3 v+c_4\extd r)\tens (v-\lambda\extd r)  \]
are
\[\und{\dim}_{\Omega^1}={2+b\lambda^2\over 1+b\lambda^2},\quad  (\id\tens (\ ,\ ))(R_{\Omega^1}\tens\id)(g)= -{1\over r^2}{\rm Vol}\left(\lambda c_1+c_2+{c_3-\lambda c_4\over 1+b\lambda^2}\right)\]
which deforms the usual trace. Note that in general in Section~3 we combined a usual trace and a cycle to give a well defined ${\rm Tr}_{\int}$ when taken together. The metric trace, by contrast, is always well-defined but in general will depend on the metric in the noncommutative case.  

(ii) We  check when we have an extendable connection. In fact the result is surprising:

\begin{lemma} Bimodule connections of the above homogeneous form are extendable if and only they are flat. In this case
\[\sigma_{\Omega^1}(\extd r\tens{\rm Vol})={\rm Vol}\tens((1+\lambda\beta)\extd r+\lambda\alpha v),\quad \sigma_{\Omega^1}(v\tens {\rm Vol})={\rm Vol}\tens(\lambda\beta'\extd r+(1+\lambda\alpha')v).\]
\end{lemma} 
\proof We set 
\[ \sigma_{\Omega^1}(\extd r\tens v\wedge\extd r)=(\wedge\tens\id)(\id\tens\sigma_{\Omega^1})(\sigma_{\Omega^1}\tens\id)(\extd r\tens v\tens\extd r)\] etc, and find the values stated, while $\sigma_{\Omega^1}((\ )\tens(\extd r)^2)=0$, $\sigma_{\Omega^1}((\ )\tens (\extd r)\wedge v)=-\sigma_{\Omega^1}((\ )\tens{\rm Vol})$ when computed in the same way. The only extendability constraints then come from $\sigma_{\Omega^1}((\ )\tens v^2)=\lambda\sigma_{\Omega^1}((\ )\tens {\rm Vol})$ which gives 4 equations
\[ (\beta-\alpha' ) \gamma + \alpha (1 - \delta +\gamma'  +  \lambda\alpha' + \lambda\beta)=0,\quad  \beta + \alpha\delta'-\gamma\beta'+\lambda\beta^2 +\lambda\alpha\beta' =0\]
\[ \alpha' - \alpha\delta' + \gamma\beta' +  \lambda\alpha'{}^2 + \lambda\alpha \beta' =0,\quad  (\alpha' - \beta) \delta' + \beta' (1 + \delta - \gamma'+ \lambda  \alpha'+  \lambda \beta)=0\]
for extendability. Remarkably, the four expressions here being set to zero are exactly the coefficients $c_1,c_2,c_3,c_4$ for the curvature listed in \cite{BegMa2}. \endproof

Using Mathematica, we find  two 4-parameter moduli of bimodule connections of the homogeneous form that are flat/extendable and non-singular in the classical limit: 

(1) We take $\alpha,\beta,\gamma,\delta$ as the parameters and 
\[ \alpha'= -\beta,\quad\beta'=  -{\beta^2\over\alpha},\quad \gamma'= \delta-1-{2 \beta \gamma\over\alpha},\quad \delta'=-{\beta (\alpha +\beta \gamma)\over\alpha^2}\]

(2) We take  $\gamma,\gamma',\delta,\delta'$ as the parameters and  $\alpha=\alpha'=\beta=\beta'=0$.

So these are all the connections of our restricted form that make $(\Omega^1,\nabla_{\Omega^1},\sigma_{\Omega^1})$ objects  of ${}_A\CG_A$ and are non-singular as $\lambda\to0$, and all necessarily flat. 

(iii)  Next we ask which of our two 4-parameter moduli of nonsingular extendable connections in (ii) are metric compatible. In fact the family (1) have no intersection with metric compatibility for generic or nonsingular $b$. The family (2) has a 1-parameter sub-family of metric-compatible connections, namely  $\alpha=\alpha'=\beta=\beta'=0$ and
\[
\gamma'=-\delta=\frac{\gamma\,\lambda}{2}\ ,\quad \delta'=-\,\frac{\gamma(1+b\,\lambda^2)}{b}
\]
where we take $\gamma$ to be the free variable. Their torsion is
always nonzero,
\[
T(\extd r)=-\gamma\,r^{-1}\,\mathrm{Vol}\ ,\quad T(v)=2\,r^{-1}\,\mathrm{Vol}
\]
and the only $*$-preserving connection in this 1-parameter family is the one with $\gamma=0$. This unique point has $\nabla_{\Omega^1}(\extd r)=\nabla_{\Omega^1}(v)=0$ and hence exactly quantises the flat metric-compatible Poisson-compatible connection which underlies the quantisation of the differential calculus and metric in Poisson-Riemannian geometry\cite{BegMa3}. This is the quantum connection $\nabla_Q$ constructed in the latter paper at semiclassical order.  In particular, the unique $*$-preserving non-singular quantum Levi-Civita connection in \cite{BegMa2} cannot be extendable by our above analysis as it is not flat. 

(iii) As is often the case in noncommutative geometry, one may need something weaker namely cotorsion-free extendable connections. Here of type (1) we have two cases

 \[{\rm (1a):}\quad 
\gamma=-\frac{\alpha  \left(2 \alpha   b  \lambda ^2+\alpha +\beta  b 
   \lambda + b \right)}{ b  (\alpha  \lambda +\beta)},\quad \delta=-\frac{\alpha  \left(2 \beta  b  \lambda ^2+\beta+2  b  \lambda \right)+\beta
    b  (\beta \lambda +3)}{ b  (\alpha  \lambda +\beta)}
   \]
where $\alpha,\beta$ are the parameters, and
\[{\rm (1b):}\quad  \beta=-\alpha\lambda,\quad \delta={\alpha\over b}  \left(1+b\lambda ^2\right)-\gamma  \lambda
   -1
   \]
where $\alpha,\gamma$ are the parameters. Of type (2) we have another two-parameter family namely $\alpha=\alpha'=\beta=\beta'=0$ and
\[
\gamma'=\gamma \lambda -2,\quad \delta'=(\delta+2) \lambda
\]
where $\gamma,\delta$ are the parameters. Thus, we have three 2-parameter moduli of extendable cotorsion free connections, all of them necessarily flat.  

If we now also impose torsion-free then we have a unique torsion free cotorsion-free or `weak quantum Levi-Civita' extendable connection of type (1b) with 
\[\alpha={4b\over 1+b\lambda^2},\quad \gamma=0.\]
In the classical limit this gives us a natural connection 
\[ \nabla \extd r={1\over r}(4b v\tens v+ 3 \extd r\tens\extd r),\quad \nabla v={2\over r}\extd r\tens v\]
on our curved manifold which is torsion free and weakly metric compatible in the sense $(\wedge\tens\id)\nabla g=0$. Unlike the Levi-Civita connection, this one is extendable on quantisation.

\subsection{Extended connections on the $q$-sphere} \label{grassspherehopf} 

We study  $A=\C_q[S^2]$ the standard $q$-sphere appearing in the base of the $q$-Hopf fibration cf\cite{BrzMa:gau,Ma:spin} as the degree 0 component of $\C_q[SL_2]$ for a $\Z$-grading corresponding to a diagonal coaction of $\C[t,t^{-1}]$. Here the standard quantum group generators  $a,c$ have grade 1 and $b,d$ grade -1 while the $q$-sphere generators are $z=cd=q^{-1}dc, z^*=-qab=-ba, x=-q^{-1}bc$ with the inherited relations
\[  zz^*=q^2x(1-q^2x),\quad z^*z=x(1-x),\quad z x=q^2 x z,\quad z^*x=q^{-1}x z^*.  \]
The quantum principal bundle here has associated quantum vector bundles or projective modules  $E_n$ for each integer $n$ realised as the grade $-n$ component of $\C_q[SL_2]$. These are each bimodules by multiplication in the quantum group and have an induced charge $n$  $q$-monopole connection $\nabla_{E_n}$  with curvature $R_{E_n}(e)= q^3[n]_{q^2}{\rm Vol}\tens e$, see \cite{Ma:spin,Ma:ltcc},  where we used $[n]_{q^2}=(1-q^{2n})/(1-q^2)$. This is a bimodule map as ${\rm Vol}$ is central. It is also possible to show that $\nabla_{E_n}$ is a bimodule connection. This was already explained for charge $\pm1$ in \cite{BMdirac} in our geometric realisation of the $q$-Dirac operator but the same argument applies also for all $n$. Then Lemma~\ref{dinnu} tells us that $\nabla_{E_n}$ is extendable with
\[ \sigma_{E_n}(f\tens{\rm Vol})=q^{2n}{\rm Vol}\tens f\]
for $f\in E_n$. In this way $(E_n,\nabla_{E_n},\sigma_{E_n})$ are objects of  ${}_{\C_q[S^2]}\mathcal{G}_{\C_q[S^2]}$.  Clearly the curvature 2-form according to Proposition~\ref{promisedphi} is $\omega_{E_n}=q^3[n]_{q^2} {\rm Vol}$.  

To give details of the braiding we will focus on  $E_1=\C_q[SL_2]_{-1}$ with module generators $b,d$ and the Grassmann connection for a standard choice of projector matrix,  which equivalently constructs the $q$-monopole. This can be written as  cf. \cite{Ma:ltcc},
\[\nabla_{E_1} f=\extd(f a)\tens d-q^{-1}\extd(fc)\tens b\]
for all $f\in E_1$. This is a bimodule connection with
\begin{align*}
\sigma_{E_1}(b\tens\extd z)=&\ \nabla_{E_1}(bz)-(\nabla_{E_1} b)z =  \extd(bza)\tens d - q^{-1}\extd(bzc)\tens b
-  \extd(ba)\tens dz + q^{-1}\extd(bc)\tens bz \\
=&\  \extd(cdab)\tens d - \extd(bc z)\tens b
- q \, (\extd(ba))z\tens d + q^{-2}(\extd(bc))z\tens b\\
=&\  -q^{-1}\extd(zz^*)\tens d + q\, \extd(xz)\tens b
+ q \, (\extd z^*)z\tens d - q^{-1}(\extd x)z\tens b \cr
=&\ (q\, \extd(xz)- q^{-1}(\extd x)z)\tens b + (q \, (\extd z^*)z-q^{-1}\extd(zz^*)) \tens d \\
=&q x\extd z\tens b+q((q^2-1)x\extd x- z^*\extd z)\tens d
\end{align*}
where we move $z$ from the right using the commutation relations  in $\C_q[SL_2]$ and through  $\tens_A$. We then identify resulting expressions in terms of the $q$-sphere generators, use the module relations 
\[ (x-1)b=z^*d,\quad z b=-q^2 x d\]
and the $q$-sphere relations for the final form. In the same way one can find
\begin{align*} \sigma_{E_1}(b\tens \extd z^*)&=-q\extd(z^*)^2\tens d+q((\extd z^*)(q^2 x-1)+q^{-2}\extd(x z^*)-(\extd x)z^*)\tens b\\ 
&=q(x\extd z^*-(q^2-1)z^*\extd x)\tens b-qz^*\extd z^*\tens d\\
\sigma_{E_1}(b\tens\extd x)&=x\extd x\tens b-q^2z^*\extd x\tens d.\end{align*}
These $\sigma_{E_1}$ are necessarily compatible with the sphere projector relation
\[qz^*\extd z+q^{-1}z\extd z^*+q^2((2)_qx-q^{-1})\extd x=0\]
where $(2)_q=q+q^{-1}$, which is a nice check of the formulae using the holomorphic calculus in \cite{Ma:spin}. This gives us the charge 1 $q$-monopole more explicitly as an object of  ${}_{\C_q[S^2]}\mathcal{G}_{\C_q[S^2]}$. The curvature 2-form according to Proposition~\ref{promisedphi} comes out as $\omega_{E_1}=q^3{\rm Vol}$. 

We also conclude that $\Omega^1\isom E_2\oplus E_{-2}$ from \cite{Ma:spin} together with the induced connection is an object of ${}_{\C_q[S^2]}\mathcal{G}_{\C_q[S^2]}$. This connection $\nabla_{\Omega^1}$ is torsion-free  and metric-compatible or `quantum Levi-Civita' for the metric
\[ g=q\extd z^*\tens\extd z+q^{-1}\extd z\tens\extd z^*+q^2(2)_q\extd x\tens\extd x\]
as essentially shown in \cite{Ma:spin}. That work did not consider $\nabla_{\Omega^1}$ as a bimodule connection and hence only showed the weaker cotorsion freeness, but full metric compatibility can be checked after constructing the braiding by the method in \cite{BMdirac}. We can therefore apply Corollary~\ref{yyuupp} and Corollary~\ref{curvtrace}. Here\cite{Ma:spin} 
\[ R_{\Omega^1}(\del f)=q^4(2)_q{\rm Vol}\tens \del f,\quad R_{\Omega^1}(\bar\del f)=-(2)_q{\rm Vol}\tens\bar\del f\] 
for $f\in \C_q[S^2]$ so that
\[ (R_{\Omega^1}\tens\id+(\sigma_{\Omega^1}\tens\id)(\id\tens R_{\Omega^1}))g=0\]
as one can also see directly given that $g$ has only mixed terms in the decomposition into holomorphic and holomorphic parts. For example, if the first tensor factor of a part of $g$ is holomorphic, we have $q^4(2)_q$ from the first term. $R_{\Omega^1}$ in the second term is $-(2)_q$ and $\sigma_{\Omega^1}$ acting on the holomorphic factor $\tens{\rm Vol}$ gives $q^4$ as this part is in $E_2$. This illustrates the antisymmetry of the Riemann tensor in Corollary~\ref{yyuupp}. There are no non-trivial Bianchi identities to illustrate here as $\Omega^3=0$. 

For the metric traces in Corollary~\ref{curvtrace}, we find similarly find
\[ \und\dim_{\Omega^1}=(2)_{q^2},\quad  (\id\tens (\ ,\ ))(R_{\Omega^1}\tens\id)g=(q^4-1)(2)_q(2)_{q^2}{\rm Vol}\]
where $1,{\rm Vol}$ are the generators of $\coH_{\rm dR}$.  We see that the metric traces $q$-deform their classical values. 

Finally, we check the square of the Dirac operator for the $q$-sphere in \cite{Ma:spin,BMdirac}. Here the spinor bundle is $\CS=E_1\oplus E_{-1}$ with  generators $f^\pm$, say, and the  Clifford action is defined as by nonzero values
\[
 f e^+\la( y f^-) = \alpha\, f y f^+\ ,\ h e^-\la (x f^+) = \beta\, h x f^-
 \]
for $f$ of degree -2, $h$ of degree 2,  $x$ of degree -1 and $y$ of degree 1, where $\alpha,\beta$ are some parameters as in \cite{BMdirac}. We have written the generators $e^\pm$ from the already-used identification of $E_{\pm 2}$ as graded components of $\C_q[SL_2]$, these being two of the basic forms in the 3D calculus on the quantum group as in \cite{Wor}. Two applications of $\la$ imply the only nonzero values
\[  f e^+ \la(h e^-\la (x f^+)) =\alpha\beta\, fh x f^+,\quad h e^-\la( f e^+\la( y f^-)) = \alpha\beta\, h f y f^-.  \]
In these terms the metric and inverse metric are 
\[g=-q^2e^+ D_1D_1'\tens D_2' D_2 e^--e^-\tilde D_1\tilde D_1'\tens \tilde D_2' \tilde D_2 e^+\]
\[  (f e^+ ,h e^-)  =- q^2 fh, \quad  (h e^- , f e^+)  =-  hf,\]
where we use the split $q$-determinants $D=d \tens a-q b\tens c=D_1\tens D_2$ (summation understood), $\tilde D=a\tens d-q^{-1}c\tens b=\tilde D_1\tens\tilde D_2$, and the primes denote independent copies. We also have  the volume form 
\[ {\rm Vol}=e^+D_1D_1'\wedge D_2'D_2e^-=e^+\wedge e^-\]
(the geometric normalisation is actually $\imath$ times this). Then  the condition in Lemma~\ref{dirsqu} 
is readily found to amount to 
\[  {\rm Vol}\la  (x f^+) = q^{-2} \kappa\, x f^+,\quad {\rm Vol}\la (y f^-) = - \kappa\, q^2\, y f^- \]
\[ \kappa={q^{-1}\over (2)_{q^2}} \alpha\beta,\quad \varphi(xf^+)=qxf^+,\quad \varphi(yf^-)=q^{-1}yf^-\]
so ${\rm Vol}\la s=q^{2|s|}s$ and $\varphi(s)=q^{-|s|}s$ on sections when viewed in $\C_q[SL_2]$ of appropriate grade $|s|$. As the connections on the spinor bundle $\CS=E_{1}\oplus E_{-1}$ are the standard associated connections as already described for $E_n$, the curvature is
\[
R_\mathcal{S} (x\,f^++y\,f^-) = q^3[1]_{q^2}{\rm Vol} \tens x\,f^+ + q^3[-1]_{q^2}{\rm Vol} \tens y\,f^-
\]
giving us
\[
\la R_\mathcal{S}= \kappa q^2\varphi^{-1}
\]
Also, as everything comes associated from the $q$-monopole connection on the quantum principal bundle one can expect, and find, that $\doublenabla(\la)=0$. Finally, the torsion is zero as the connection on $\Omega^1$ is the quantum-Levi-Civita one. Hence Lemma~\ref{dirsqu} tells us that
\[ \dirac^2=\kappa \varphi^{-1}\circ\left( (\ ,\ )\nabla_{\Omega^1\tens\mathcal{S}}\nabla_\mathcal{S}+ q^2\varphi^{-1}\right)\]
as our $q$-Lichnerowicz formula. The last term reflects a constant curvature with an additional $\varphi^{-1}$ automorphism. One has similar results for the Dirac operator on $M_2(\C)$ in \cite{BMdirac} but with zero curvature. For the $q$-disk Dirac operator also in \cite{BMdirac}, we have a vanishing curvature contribution but nonzero $\doublenabla(\la)$, which does not preclude the possibility of a different $\la$ for this model to make it more like the $q$-sphere.

\end{document}